\documentclass[nopreprintline,12pt]{elsarticle}
\pdfoutput=1

\usepackage[percent]{overpic}
\usepackage{graphicx,amssymb,mathrsfs,graphics,algorithm,algpseudocode,indentfirst}
\usepackage{amsmath,amssymb,enumerate,color,multirow,bm}
\usepackage{graphicx}
\usepackage{epstopdf}
\usepackage{epsfig}

\newtheorem{example}{{Example}}[section]

\newcommand{\ba}{\begin{array}}
	\newcommand{\ea}{\end{array}}
\newcommand{\bt}{\begin{tabular}}
	\newcommand{\et}{\end{tabular}}
\newcommand{\btb}{\begin{table}}
	\newcommand{\etb}{\end{table}}
\newcommand{\bc}{\begin{center}}
	\newcommand{\ec}{\end{center}}
\newcommand{\bea}{\begin{eqnarray}}
\newcommand{\eea}{\end{eqnarray}}
\newcommand{\Bea}{\begin{eqnarray*}}
	\newcommand{\Eea}{\end{eqnarray*}}
\newcommand{\beq}{\begin{equation}}
\newcommand{\eeq}{\end{equation}}

\def\bigdot{\boldsymbol{\,\cdot\,}}
\def\eref#1{{\rm (\ref{#1})}}

\def\df{\varepsilon}

\newcommand{\ve}[1]{\bm{#1}}
\def\R{\mathbb{R}} 
\def\M{\mathcal{S}}
\def\Hsp{\mathcal{H}}

\def\call{\mathcal{L}} 

\def\Lhs{{\call}}
\def\Rhs{g}

\def\ven{\ve{{n}}}
\def\vex{\ve{x}}
\def\vev{\ve{v}}

\def\cala{{\mathcal{A}}}

\def\qed{~\relax\ifmmode\hskip2em \Box
	\else\unskip\nobreak\hskip1em \hfill$\Box$
	\fi \newline}

\def\bigdot{\boldsymbol{\,\cdot\,}}
\def\eref#1{{\rm (\ref{#1})}}

\def\df{\varepsilon}

\def\R{\mathbb{R}} 
\def\M{\mathcal{S}}
\def\Hsp{\mathcal{H}}

\def\call{\mathcal{L}} 

\def\Lhs{{\call}}
\def\Rhs{g}

\def\ven{\ve{{n}}}
\def\vex{\ve{x}}
\def\vev{\ve{v}}

\def\cala{{\mathcal{A}}}

\def\qed{~\relax\ifmmode\hskip2em \Box
	\else\unskip\nobreak\hskip1em \hfill$\Box$
	\fi \newline}

\makeatletter
\newenvironment{breakablealgorithm}
{
	\begin{center}
		\refstepcounter{algorithm}
		\hrule height.8pt depth0pt \kern2pt
		\renewcommand{\caption}[2][\relax]{
			{\raggedright\textbf{\ALG@name~\thealgorithm} ##2\par}%
			\ifx\relax##1\relax 
			\addcontentsline{loa}{algorithm}{\protect\numberline{\thealgorithm}##2}%
			\else 
			\addcontentsline{loa}{algorithm}{\protect\numberline{\thealgorithm}##1}%
			\fi
			\kern2pt\hrule\kern2pt
		}
	}{
		\kern2pt\hrule\relax
	\end{center}
}
\makeatother

\begin{document}

\title{
Kernel-based collocation methods for heat transport on evolving surfaces 
}\tnotetext[mytitlenote]{
This work was supported by a Hong Kong Research Grant Council GRF Grant.
}
	
	\author[mymainaddress]{Meng Chen
}
	\ead{chenmeng@astri.org}
	
	\author{Leevan Ling\fnref{my2}}
	\ead{lling@hkbu.edu.hk}

	\address[mymainaddress]{Hong Kong Applied Science and Technology Research Institute Company Limited}
	\address[my2]{Department of Mathematics, Hong Kong Baptist University, Kowloon Tong, Hong Kong}

	\date{\today}
	\begin{abstract}
We propose algorithms for solving convective-diffusion partial differential equations (PDEs), which model surfactant concentration and heat transport on evolving surfaces, based on intrinsic kernel-based meshless collocation methods. The algorithms can be classified into two categories: one collocates PDEs directly and analytically, and the other approximates surface differential operators by meshless pseudospectral approaches. The former is specifically designed to handle PDEs on evolving surfaces defined by parametric equations, and the latter works on surface evolutions based on point clouds. After some convergence studies and comparisons, we demonstrate that the proposed method can solve challenging PDEs posed on surfaces with high curvatures with discontinuous initial conditions with correct physics.
	\end{abstract}
	
	
\maketitle	

	\section{Introduction}
	
	\noindent
Heat transport
on moving surfaces are modelled by convection-diffusion
	partial differential equations (PDEs) that
frequently appear in physical and biological fields.
Applications include
dealloying by surface dissolution \cite{eilks2008numerical}, pattern formation on evolving biological surfaces \cite{barreira2011surface}, modelling geometric biomembranes \cite{elliott2010modeling}, and cell motility and chemotaxis~\cite{elliott2012modelling}. The formal problem statement is as follows.
Let $\M(t)\subset \R^d$ be a continuously evolving surface with codimension one; for each $t\in[0,T]$, we also assume that the surface is smooth, closed, connected, and complete.
Furthermore, let $\ven=\ven(\bigdot,t)$ and  $\ve{\tau}=\ve{\tau}(\bigdot,t)$ represent the unit normal  and the tangent vector of $\M(t)$.

We consider convective-diffusion equations defined on $\M(t)$ in the form of
\begin{equation}\label{eq_CDpupt}
	\begin{array}{ll}
	\partial_t  u +\vev\bigdot\nabla u+(\nabla_\M \bigdot\vev) u -\df \Delta_\M u=f &\mbox{\quad in }\M(t)\times[0,T]
	\end{array}
\end{equation}
subject to  homogeneous Dirichlet initial conditions
\[
u(\bigdot,0)=u_0  \mbox{\quad on }  \M(0),
\]
where
$\df\geq0$ is the  diffusive coefficient, and
$\vev:=\vev_{n} + \vev_{\tau}$ is the velocity of surface motion whose  normal and tangent components
are denoted by
$\vev_{n}$ and $\vev_{\tau}$  respectively.
The surface differential operators in \eref{eq_CDpupt} are defined as
\[
    \nabla_\M :=(I_d-\ven\ven^T)\nabla
\]
and 
\[
    \Delta_\M:=\nabla_\M\bigdot \nabla_\M
\]
where $I_d$ is the $d$-dimensional identity matrix and $\nabla$ is the standard gradient operator in Euclidian spaces.

The problem we considered includes several physical models in the literature.
In the absence of  diffusive flux and source, i.e., $\df=0$ and $f=0$, the viscosity solution $u$ to \eref{eq_CDpupt} as $\df\rightarrow 0^+$ satisfies a mass conservation law  \cite{stone1990simple}, namely
\begin{equation}\label{eq_conlaw1}
		\frac{\text{d}}{\text{d} t}\displaystyle\int_{\M(t)}u\ \text{d}\M=0.
\end{equation}
Readers can find more details in 
\cite{barreira2011surface,dziuk2007finite,dziuk2010eulerian,grande2014eulerian}.


To solve the target PDE \eref{eq_CDpupt} on the evolving surface, a surface finite element method (SFEM) were proposed with the use of weak and variational formulas  based on  triangular meshes  intrinsically defined on surfaces \cite{dziuk2007finite,intrinsicGerhard}. Embedding techniques for static problems can also be  extended to deal with surface PDEs posed in higher-dimensional  \cite{closetpointmethod,narrowbM,Cheung+Ling-Kernembemethconv:16,marz2012calculus,Marz}. Some domain-type time-dependent FEMs have been implemented in narrow bands containing moving surfaces \cite{dziuk2010eulerian,grande2014eulerian,deckelnick2014unfitted,elliott2010numerical}. Implicit time-stepping schemes have also been used for temporal discretization in the FEM literature mentioned above, except in \cite{grande2014eulerian} where weak forms were utilized in both space and time.
In this paper, we apply kernel-based meshless collocation techniques to solve \eref{eq_CDpupt}. The proposed method is intrinsic, which means that no embedding is required. One variant of our method naturally works with point clouds, which makes it capable to handle solution-driven surface evolutions.

The rest of this paper is organized as follows. Firstly, we deal with the case where the evolutions of surfaces are prescribed by some parameterizations that are known a prior.
Section~\ref{sect_Conveq_parS} contains a novel algorithm for this case. Next, we focus on cases where the surfaces are defined by some point clouds.
In Section~\ref{AKMtheta}, we propose a second algorithm that does not require analytic information about normals and velocities of surface motions.
In Section~\ref{SecNumexp},
we first test accuracy and convergence of our first algorithm. Next, we test the proposed algorithms against mass balance in the case of \eref{eq_conlaw1}.
Lastly, as a pilot study, we examine the robustness of our method in dealing with merging surfaces, which involves with high curvatures.
We conclude the results and observations in Section \ref{sec;colclusion}.
	
	\section{A collocation method for PDEs posed on parameterized surfaces with prescribed evolution}\label{sect_Conveq_parS}
	
To begin, we consider the case that the surface $\M(t)$ and its motion are both given by some time-dependent parameterization.
Assume the evolving surface is defined as
\begin{equation}\label{eq_X}
	\M(t):=\Big\{\vex (\ve{\varphi} ,t)\in\R^d, \ \text{for} \;t\in[0,T],\ \ve{\varphi}\in \mathcal{B}\subset \R^{d-1}\Big\},
\end{equation}
where $\ve\varphi\in\mathcal{B}$ are parameters in certain parameter space that define surface points $\vex(\ve{\varphi},t)=\big[x_i(\ve{\varphi},t)\big]_{i=1}^d$ for any time $t\in[0,T]$.
The following example aims to clarify our notations.

\subsubsection*{Example:}
Consider an evolving circle with changing radius. Using  polar coordinate, we can pick $\varphi=\ve\varphi\in\mathcal{B}=[0,2\pi]\subset\R$ and for some function $r(t)>0$ to define surface points on $\M(t)$ at any time $t\geq0$ by
\[
    \vex(\varphi,t) = [x_1(\varphi,t),x_2(\varphi,t)]= [ r(t) \cos\varphi,\,r(t) \sin\varphi] \in \M(t)\subset\R^2.
\]
Obviously, the velocity for this surface motion is given by the partial derivative of $\vex$ with respect to $t$.
\qed

In our notations, we have velocity vector in \eref{eq_CDpupt} given by
\begin{equation}\label{eq_vel}
	\vev(\ve{\varphi} ,t)=\frac{\partial }{\partial t}  \vex(\ve{\varphi},t) =\begin{bmatrix}\partial_t x_i(\ve{\varphi},t)\end{bmatrix}_{i=1}^d.
\end{equation}
With all required analytic information about the evolving surface ready,  we can start discretizing \eref{eq_CDpupt} first in time and then in space.

Firstly, we semi-discretize \eref{eq_CDpupt} based on some  partition $\{t_m\}_{m=0}^N$ of  the  interval $[0,T]$. For simplicity, we assume the partition is equispaced with time stepping size $\triangle t$ and we will employ the standard $\theta$-scheme. Note that the problem \eref{eq_CDpupt} in hand  is given in the form of ${\partial_t u}=F(u)$, which is \emph{not} yet ready for applying the $\theta$-scheme.
Instead, we implicitly discretize the material derivatives  \cite{dziuk2010eulerian,deckelnick2014unfitted} of the solution to \eref{eq_CDpupt} as follows.

%

We abbreviate $\M(t^m)=:\M^m$ the snapshot of our evolving surfaces at time $t_m$ for $m=1,\ldots,T/\triangle t$, on which we have surface points
\[
    \vex^m=\vex(\ve{\varphi},t^m)\subset \M^m, \quad \ve\varphi\in\mathcal{B}.
\]
Let us also denote the solution $u$ and right-hand function $f$ at $t_m$ as
\[
u^m :=	u (\vex^m,t^m)=u\Big(\vex(\ve{\varphi},t^m),t^m\Big)
\mbox{\quad and \quad}
f^m:=f(t_m).
\]
We  approximate the material derivative of $u$ and  obtain
\begin{eqnarray*}\label{eq_dudtap}
\frac{\partial}{\partial t} u 
+  \vev \bigdot \nabla u
&=&
\frac{\text{d}}{\text{d}t} u\Big(\vex(\ve{\varphi} ,t),t\Big)
\\
&\approx&	\frac{u(\vex^{m},t^{m})-u(\vex^{m-1},t^{m-1})}{\triangle t}		
=\frac{u^{m}-u^{m-1}}{\triangle t}.		
\end{eqnarray*}
Thus, the time discretization of \eref{eq_CDpupt} by $\theta$-method is given as
	\begin{eqnarray}\label{eq_theta0}
	\frac{u^{m}-u^{m-1}}{\triangle t} +(1-\theta) \cala^m u^{m}  +\theta \cala^{m-1} u^{m-1}
	=(1-\theta)f^{m}+\theta f^{m-1},
	\end{eqnarray}
where the surface elliptic operator $\cala^m$ is defined as
	\begin{equation}\label{eq_lst}
	\cala^m :=\vev\bigdot\nabla_{\M^m} - \df \Delta_{\M^m}
	\end{equation}
on the surface $\M^m$ and the velocity $\vev=\vev(\bigdot,t^m)$ is to be evaluated exactly.
We can further simplify \eref{eq_theta0} to isolate the unknown solution $u^{m}$ and get a time-independent second-order surface PDE:
\begin{equation}\label{eq_theta}
	\Lhs^{m} u^{m} := \Big(1+{(1-\theta)}\triangle t\cala^{m}  \Big) u^{m}= \Rhs^{m},
\end{equation}
where the  right-hand function $\Rhs^m$, which  depends  on the known functions $u^{m-1}$ and $f$, is known and given by
\begin{equation}\label{hm}
\Rhs^{m}:= \Big(1-{\theta}\triangle t\cala^{m-1}\Big)u^{m-1}+\triangle t\Big((1-\theta)f^m+\theta f^{m-1}\Big).
\end{equation}
Aiming towards fully discretized problems, we employ our previously proposed overdetermined Kansa method \cite{Chen-IntrmeshmethPDEs:18} to solve surface PDEs.
We set up an overdetermined kernel-based collocation method \cite{cls} to solve surface PDE  $\Lhs^{m} u^{m} = \Rhs^m$ on surfaces
$\M^{m}$ and $\M^{m-1}$ (due to {$\Lhs^m$ and}  $\Rhs^m$)
for solution $u^{m}$.

Based on the theories in \cite{Fuselier}, we can take some commonly used symmetric positive definite (SPD) kernel and restrict it to any one of the surfaces $\M^m$ to obtain an SPD surface kernel
\[
    \Psi(\bigdot,\bigdot):\M^m\times \M^m\rightarrow \R
\]
with a certain smoothness order $\mu>(d-1)/2$.
Let $\varphi_Z\subset\mathcal{B}$ be a set of $n_Z$ points in the parameter space; similarly, let $\varphi_X\subset\mathcal{B}$ be another set of $n_X > n_Z$ points.
For each $m$, let $Z^m=\vex(\varphi_Z,t^m)$ and $X^m=\vex(\varphi_X,t^m)$ be the set of trial centers and collocation points.
We shall use denser collocation points because of available convergence theories for collocation methods in domains \cite{cls} and on surfaces \cite{Cheung+Ling-Kernembemethconv:16}. In both works, they require $X$ to be sufficiently denser with respect to $Z$ {resulting in \textit{overdetermined formulas}}, in order to establish stability estimates\footnote{In this context, stability means that errors can be bounded in continuous norms with discrete residuals.}.

Trial functions at time $t^m$ are expressed by a linear combination
\begin{equation}\label{eq_intS}
	u^m = \Psi(\bigdot,Z^m)\boldsymbol{\lambda}_Z^m
	:= \sum_{z_i\in Z^m}\lambda_i^m \Psi(\bigdot,z_i),
\end{equation}
where $\boldsymbol\lambda_Z^m:=\{\lambda_i^m\}_{i=1}^{n_Z}$ is the unknown coefficient vector to be determined. As in the original Kansa method, we collocate \eref{eq_theta} and \eref{hm} at set $X^m$ to yield an $n_X\times n_Z$ matrix system
\begin{equation}\label{kansa}
    [\call^m \Psi(X^m,Z^m)]\boldsymbol\lambda_Z^m = \Rhs^m_{|X}.
\end{equation}
In the Euclidean counterpart, overdetermined formulations guarantee stability {as mentioned before}
, which ultimately leads to error estimates \cite{Cheung+Ling-Kernembemethconv:16,cls}.

Given the parametric equation \eref{eq_X}, the surface differential operator in \eref{eq_lst} can be rewritten without any implicit dependency on the surface.
Thus, $\call^m \Psi:\M^m\times \M^m\to \R$ can also be analytically evaluated for each $m=1,\ldots,T/\triangle t$.
The right hand vector $\Rhs^m_{|X}=\Rhs^m(X^m)$ is the nodal value of \eref{hm} evaluated at $X^m$.
The overdetermined system \eref{kansa} should be solved in the least-squares sense.
	We summarize with Algorithm~\ref{KPM-CN}.
	
	\begin{breakablealgorithm}
		\caption{: For solving PDE \eref{eq_CDpupt} posed on a parametrized surface \eref{eq_X}}
		\label{KPM-CN}
		\begin{algorithmic}
			\State{${\bf Initialization:}$}
			\State{$\qquad set ~ t=t^0,\ select~ 0\leq \theta\leq 1;$}
            \State{$\qquad define ~\varphi_Z, \varphi_X\subset \mathcal{B},\ compute~Z^0, X^0\subset \M^0;$}
            \State{$\qquad interpolate\ u_0\ at\ Z^0\ to\ obtain\ \boldsymbol\lambda_Z^0;$}
			\State{$\qquad symbolically~derive~analytic ~formulas~of~\cala^m~in~\eref{eq_lst};$}
            \State{$\qquad m=1;$}

			\State{${\bf While}~(t\leq T)$}
			\State{$\qquad compute~\{ X^{m},Z^{m}\}\subset\M^{m};$}
			\State{$\qquad assemble~and\ solve~\eref{kansa}\ in\ the\ least$-$squares\ sense;$}
			\State{$\qquad m=m+1;$}
			\State{$\qquad t=t^0+(m+1)\triangle t;$}
			\State{${\bf End}$}
		\end{algorithmic}
	\end{breakablealgorithm}


\section{Another collocation method for point clouds}\label{AKMtheta}

\begin{figure}
		\centering
		\begin{overpic}[width=0.46\textwidth,trim=85 40 85 25, clip=true,tics=10]{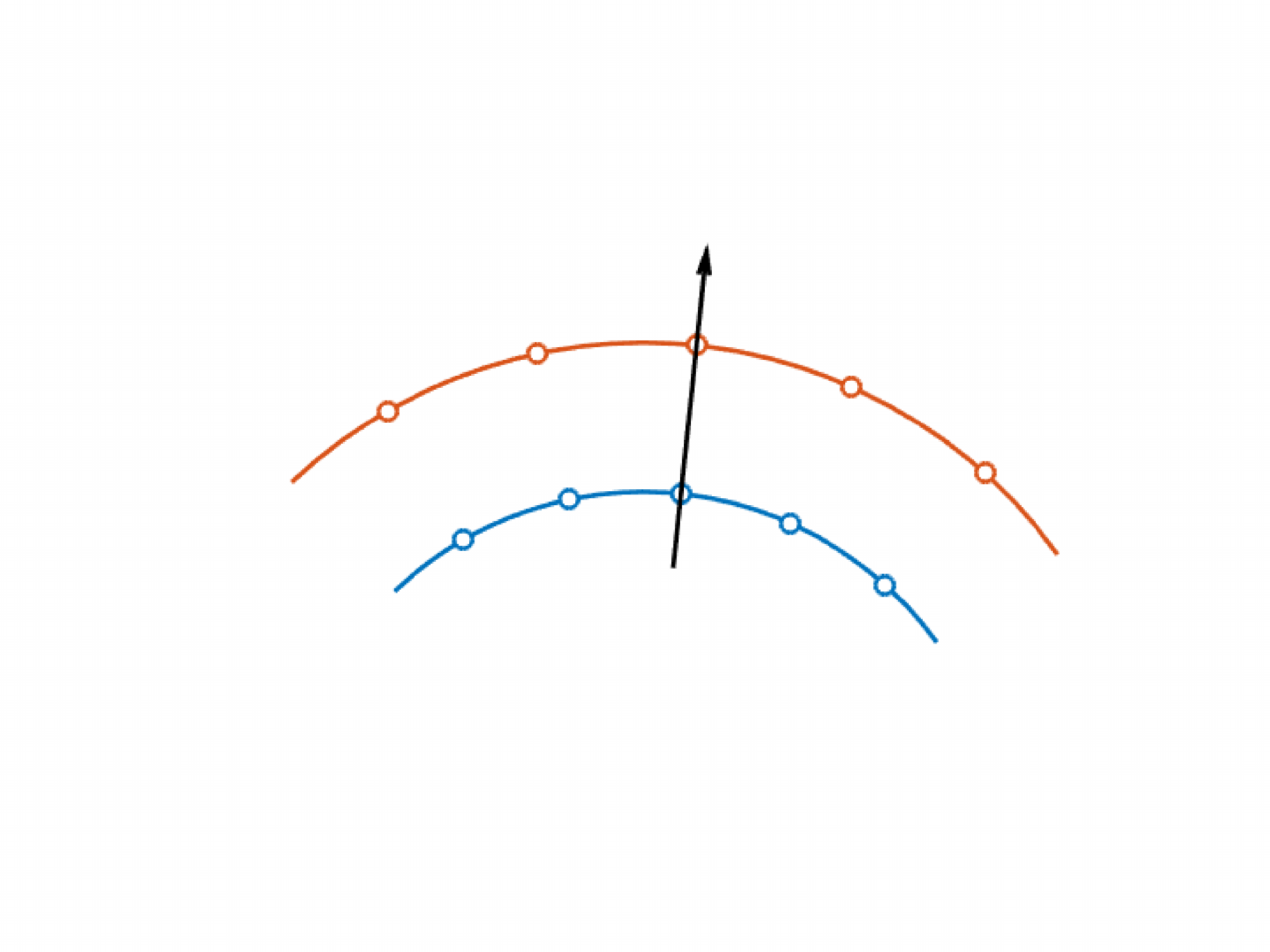}
			
			\put(59,61) {\scriptsize  ${\vex}^{m+1}$}
			\put(56,36) {\scriptsize  $\vex^m$}
			\put(80,20) {\scriptsize  $\vex(\ve{\varphi},t )$}
			\put(48,47) {\scriptsize  $\vev$}
			\put(8,38) {\scriptsize  $t^{m+1}$}
			\put(20,25) {\scriptsize  $t^{m}$}
		\end{overpic}
		\caption{The point-to-point motion of a point $\vex(\ve{\varphi},t)$ on a surface from $t^m$ to $t^{m+1}$  with the velocity  $\vev$.}\label{Fig_dudt}
\end{figure}
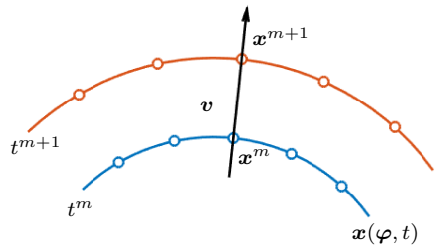

\begin{figure}
		\centering
		
		\begin{overpic}
			[width=0.48\textwidth,trim=110 80 115 65, clip=true,tics=10]{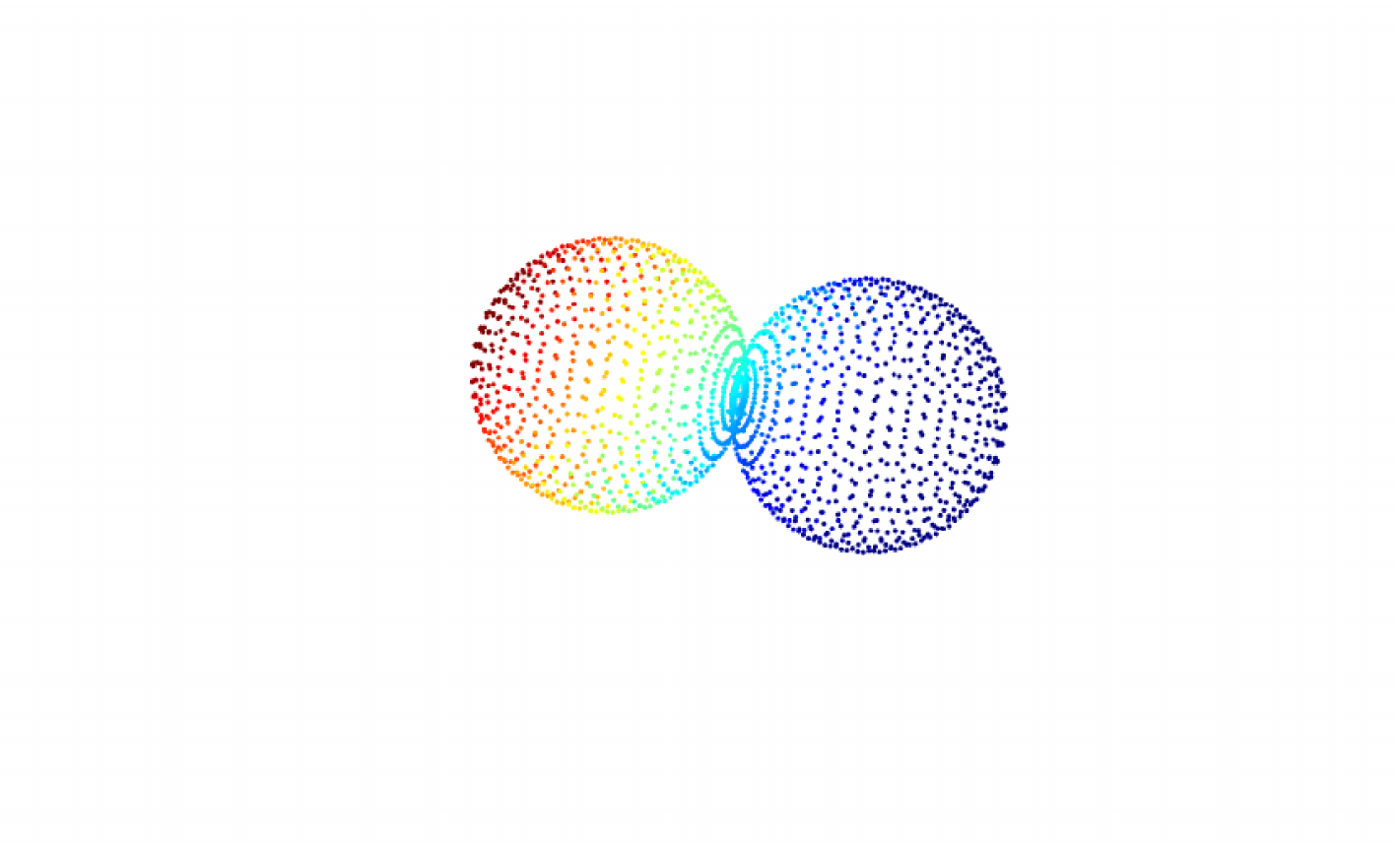}
			\put(45,-3) {\scriptsize  (a) $t=0$}
		\end{overpic}
		\begin{overpic}
			[width=0.48\textwidth,trim=110 80 115 65, clip=true,tics=10]{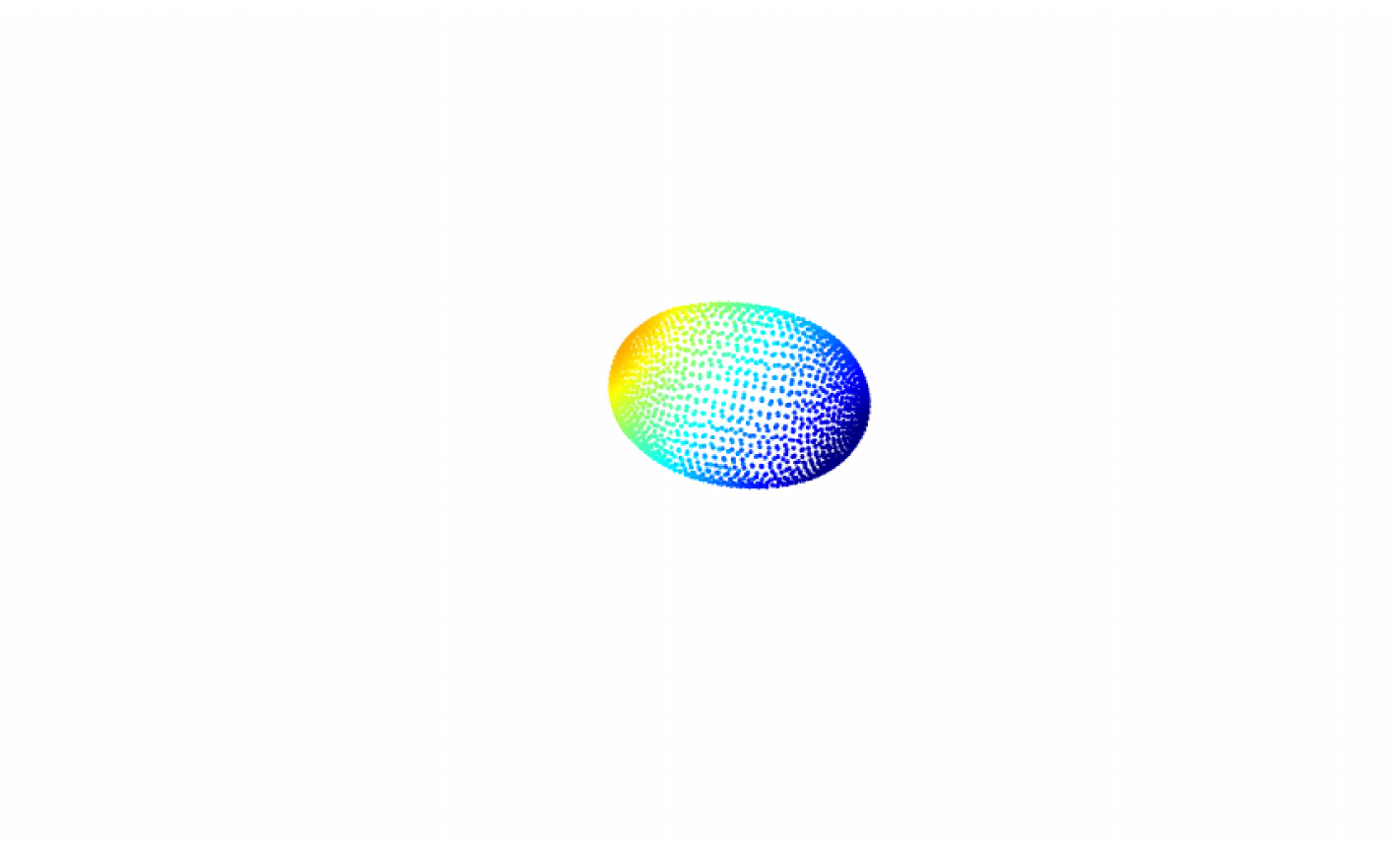}
			\put(45,-3) {\scriptsize  (b) $t=1.5$}
		\end{overpic}
		
		\caption{The snapshots of an evolving surface based on point clouds by mean-curvature motion.}\label{Fig_MovMeanCur}
\end{figure}

In the case of surfaces defined by point clouds, some approximated method is required to discretize space without analytic formulas for the normal vectors. In this section, we assume that the initial surface $\M^0$ is given by a set of points in $\R^d$.
We also assume the surface evolution velocity $\vev=\vev(\vex,t)$ is independent to the PDE solution.
As surface points $\vex(\ve{\varphi},t)$ evolving from $\M^m$ to $\M^{m+1}$  ($m\geq0$) with velocity  $\vev$, we can approximate new  points ${\vex}^{m+1}$ on $\M^{m+1}$ by
\begin{equation}\label{eq_ptp}
		{\vex}^{m+1}\approx{\vex}^{m}+\vev \triangle t,
\end{equation}
as shown in Figure~\ref{Fig_dudt}. This process allows us to track all surfaces and place data points $X$ and $Z$ via local interpolation on them as in Algorithm~\ref{KPM-CN}; see Figure~\ref{Fig_MovMeanCur} for a schematic demonstration.
The problem here is that  $\cala^m$ in \eref{eq_lst} can no longer be obtained analytically.
Our solution is to take the RBF pseudospectral approach \cite{Spere2}  to approximate the surface   Laplacian.
The ultimate goal is to obtain an approximated version of overdetermined linear system as in \eref{kansa}.

Recall that $\nabla_\M :=(I_d-\ven\ven^T)\nabla$. We can approximate $\ven=\ven(\vex)$ for $\vex\in\M^m$ from the point cloud $\{\vex^m_k\}_{k=1}^K$ specifying $\M^m$ by
\begin{equation}\label{find n}
  \ven= \frac{\nabla\Psi(\bigdot,Z)\Psi(Z,Z)^{-1} \mathbf{1}_{n_Z}}{\|\nabla\Psi(\bigdot,Z)\Psi(Z,Z)^{-1} \mathbf{1}_{n_Z}\|}\,:\,\M^m\to\R^d, \quad Z\subset\M^m
\end{equation}
where $\mathbf{1}_{n_Z}$ is the all-ones vector of size $n_Z$, see \cite{fornberg2002observations,marchandise2012cad}.
When implementing \eref{find n}, we use  subsets of $Z$ surrounding the points of interest for the sake of efficiency.
Once we have all normal vectors at $X$ approximated, the computation of entries in the approximated gradient matrix $\widetilde\nabla_\M(X,Z)$ is straightforward. The approximation error here is due to error in $\ven$ and the exactness of data point locations on $\M^m$.

The same technique is undesirable for the approximation of surface Laplacian $\Delta_\M:=\nabla_\M\bigdot \nabla_\M$, as the surface divergence operator will act upon $\ven$ and we simply want to avoid approximating second derivatives from point cloud.
In the RBF pseudospectral approach,
we rewrite  \eref{eq_intS} in terms of
the unknown nodal data $u_{|Z}^{{m}} := {u^m}(Z)$ of $u^m$ at trial centers $Z\subset\M^m$.
This yields
\begin{equation}\label{eq_intS2}
	u^m = \Psi(\bigdot,Z)\Psi(Z,Z)^{-1}u_{|Z}^m.
\end{equation}
We can apply the approximation and evaluate the surface gradient at $Z$ while keeping $u_{|Z}^m$ unknown to yield the following  $n_Z\times d$ matrix
\[
   (\widetilde\nabla_\M u^m)_{|Z} = \widetilde\nabla_\M\Psi(Z,Z)\Psi(Z,Z)^{-1}u_{|Z}^m.
\]
At this point, we have approximated nodal values of each component of $\nabla_\M$ at $Z$, which still depends on the unknown $u_{|Z}^m$. Without explicitly writing out the index of these components for simplicity,  we can express them similar to \eref{eq_intS2} by its interpolant
\[
   \widetilde\nabla_\M u^m =  \Psi(\bigdot,Z)\Psi(Z,Z)^{-1} (\widetilde\nabla_\M u^m)_{|Z},
\]
to which we can apply an approximate surface divergence at $X\subset \M^m$ to yield
\[
   \widetilde\Delta_\M u^m = \big(\widetilde\nabla_\M \bigdot  \Psi\big)(X,Z)\Psi(Z,Z)^{-1} (\widetilde\nabla_\M u^m)_{|Z}.
\]
This apparently complicated procedure did come with some neat simplifications \cite{fuselier2013high}. Let $\vex\in\M^m$ be a surface point on which the (approximated) normal vector is denoted by $\ven(\vex)$. Also let
\[
P(\vex) = I_d -\ven(\vex) \ven(\vex)^T =: [p_1(\vex), \ldots,p_d(\vex)] \in \R^{d\times d}
\]
be the (approximated) projection operator to $\vex\in\M^m$ with columns $p_k$, $k=1,\ldots,d$.
For each $k$, we define row-vector functions
\[
    G_k(\vex,Z) := p_k(\vex)^T [\nabla\Psi(\vex,Z)] [\Psi(Z,Z)]^{-1}.
\]
Then, we have
\begin{eqnarray}
 [ \widetilde{\nabla}_\M u^m ]_k &=&  G_k(\bigdot,Z) u_{|Z}^m , \quad k=1,\ldots,d,\label{LL1}\\
  \widetilde{\Delta}_\M u^m  &=& \sum_{k=1}^{d}G_k(\bigdot,Z)G_k(Z,Z)u_{|Z}^m, \label{LL2a}
\end{eqnarray}
or, alternatively, we can use $(\widetilde\nabla_\M u^m)_{|X}$ as data to approximate surface Laplacian by
\begin{eqnarray}\label{LL2b}
\widetilde{\Delta}_\M u^m  &=& \sum_{k=1}^{d}G_k(\bigdot,X)G_k(X,Z)u_{|Z}^m.
\end{eqnarray}
These equations, along with \eref{eq_intS2}, are sufficient for obtaining an approximated version of the overdetermined linear system in \eref{kansa} for solving
\eref{eq_theta} in terms of unknown $u_{|Z}$. Using the right most inverse of the interpolation matrix of $\Psi$ at $Z$ and the relationship $\boldsymbol\lambda_Z^m = [\Psi(Z,Z)]^{-1} u_{|Z}^m$, we can recast the system with unknown $\boldsymbol\lambda_Z^m$.
Below is a pseudocode that summarizes the  algorithm above:
Algorithm \ref{AKM-CN}, specifically denoted as Algorithm \ref{AKM-CN}a by \eref{LL2a} and Algorithm \ref{AKM-CN}b by \eref{LL2b}.

\begin{breakablealgorithm}
		\caption{a/b: For solving PDE \eref{eq_CDpupt} posed on point cloud}
		\label{AKM-CN}
		\begin{algorithmic}
			\State{${\bf Initialization:}$}
			\State{$\qquad set ~ t=t^0,\ select~ 0\leq \theta\leq 1;$}
            \State{$\qquad define~~Z^0, X^0\subset \M^0\ based\ on\ some\ point\ cloud$}
            \State{$\qquad compute~\ven(X^0)\ using\ \eref{find n};$}
            \State{$\qquad interpolate\ u_0\ at\ Z^0\ to\ obtain\ \boldsymbol\lambda_Z^0;$}
            \State{$\qquad m=1;$}

			\State{${\bf While}~(t\leq T)$}
\State{$\qquad compute~\vev~and~update~\{ X^{m},Z^{m}\}\subset\M^{m};$}			
\State{$\qquad approximate~entries~in~\eref{kansa}~using\  \eref{eq_intS2}-\eref{LL1},~and~\eref{LL2a}/\eref{LL2b};$}
			\State{$\qquad solve~\eref{kansa}\ in\ the\ least$-$squares\ sense;$}
			\State{$\qquad m=m+1;$}
			\State{$\qquad t=t^0+(m+1)\triangle t;$}
			\State{${\bf End}$}
		\end{algorithmic}
	\end{breakablealgorithm}

\section{Numerical examples}\label{SecNumexp}
We provide five  examples to numerically study the proposed algorithms for solving convective-diffusion equations  \eref{eq_CDpupt} on evolving surfaces.
In all examples, we use the restricted Whittle-Mat\'{e}rn-Sobolev kernels and  $\theta=0.5$ in \eref{eq_theta} for the Crank-Nicolson method (CN).

To   quantify the
accuracy and convergence of our numerical approximations, we employ two error measures as follows. Let $u^*$ denote the exact solution. The $L^\infty\big([0,T];L^2(\M(t))\big)$ error is defined by
	\[
	\sup_{[0,T]}|| u^* -u||_{L^2(\M(t))}.
	\]
	Similarly, the $L^\infty\big([0,T];\Hsp^1(\M(t))\big)$ semi-norm error is given by
	\[
	||\nabla_{\M} u^* -\nabla_{\M}u||_{L^\infty\big([0,T];L^2(\M(t))\big)},
	\]
	and lastly, the $L^\infty\big([0,T];\Hsp^2(\M(t))\big)$ semi-norm error is defined by
	\[
	||\Delta_{\M} u^* -\Delta_{\M}u||_{L^\infty\big([0,T];L^2(\M(t))\big)}.
	\]
At any given time $t^m\in[0,T]$, $m\geq1$, by using one of the above  errors $e^{m}$  and $e^{m-1}$ at two consecutive time steps corresponding to fill distances $h_{X^{m-1}}$ and $h_{X^{m}}$ for collocation point sets  $X^{m-1}\subset\M^{m-1}$  and $X^m\subset\M^{m}$ respectively, we  estimate the order of convergence by
	\begin{equation}\label{eoc}
	\text{eoc}^{m}:=\frac{\log(e^{m}/e^{m-1})}{\log(h_{X^{m}}/h_{X^{m-1}})}.
	\end{equation}

In cases of smooth initial conditions (ICs) on surfaces without (extremely) high curvatures, we can alternatively employ overdetermined expressions with $n_X>n_Z$ or exactly determined formulas with the same point sets ($X=Z$) for both  Algorithms~\ref{KPM-CN}  and \ref{AKM-CN}. Having that said, oversampling in Algorithm~\ref{KPM-CN} can still result in better accuracy; see Example~\ref{Ex_MovCur2}.
The situations are quite different when imposing discontinuous ICs on surfaces with high curvatures, in which oversampling is essential.
In our concluding simulations in Examples~\ref{Ex_ContSurf} and \ref{Ex_DiscSurf}, we will demonstrate the robustness of our proposed methods.


	\begin{example}\label{Ex_MovCur2}{Evolving curve in $\R^2$.}\end{example}
	In the first example, we compare Algorithm~\ref{KPM-CN} with two time-dependent FEMs \cite{grande2014eulerian,deckelnick2014unfitted} for solving \eref{eq_CDpupt} on two moving surfaces.
	We use the same diffusivity $\df=1$ and the source terms $f(\vex ,t)$ by the symbolic mean by providing
the same exact solutions $u^*(\vex,t)$ in the FEM papers. The smoothness order  is set to be $\mu=4$.

In \cite{deckelnick2014unfitted},  weak formulations of \eref{eq_CDpupt} were  discretized in some narrow bands containing the evolving surface, and an implicit temporal scheme was used to solve \eref{eq_CDpupt}, which yields an \emph{unfitted finite element method}. We aim to examine the same experiment as in \cite[Section~6.3]{deckelnick2014unfitted} with the use of Algorithm~\ref{KPM-CN}, imposed on the following evolving curve,
			\[
			 \M(\vex,t) =\frac{x_1^2}{1+0.25\sin(2\pi t)}+x_2^2-1=0\ \mbox{\quad for } t\in[0,0.5].
			\]
			
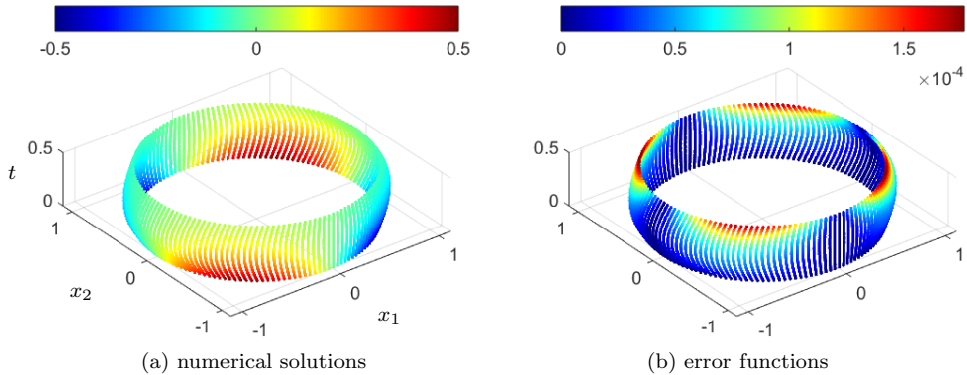
\begin{figure}
				\centering
				\begin{tabular}{cc}
					\begin{overpic}[width=0.46\textwidth,trim=85 40 85 25, clip=true,tics=10]{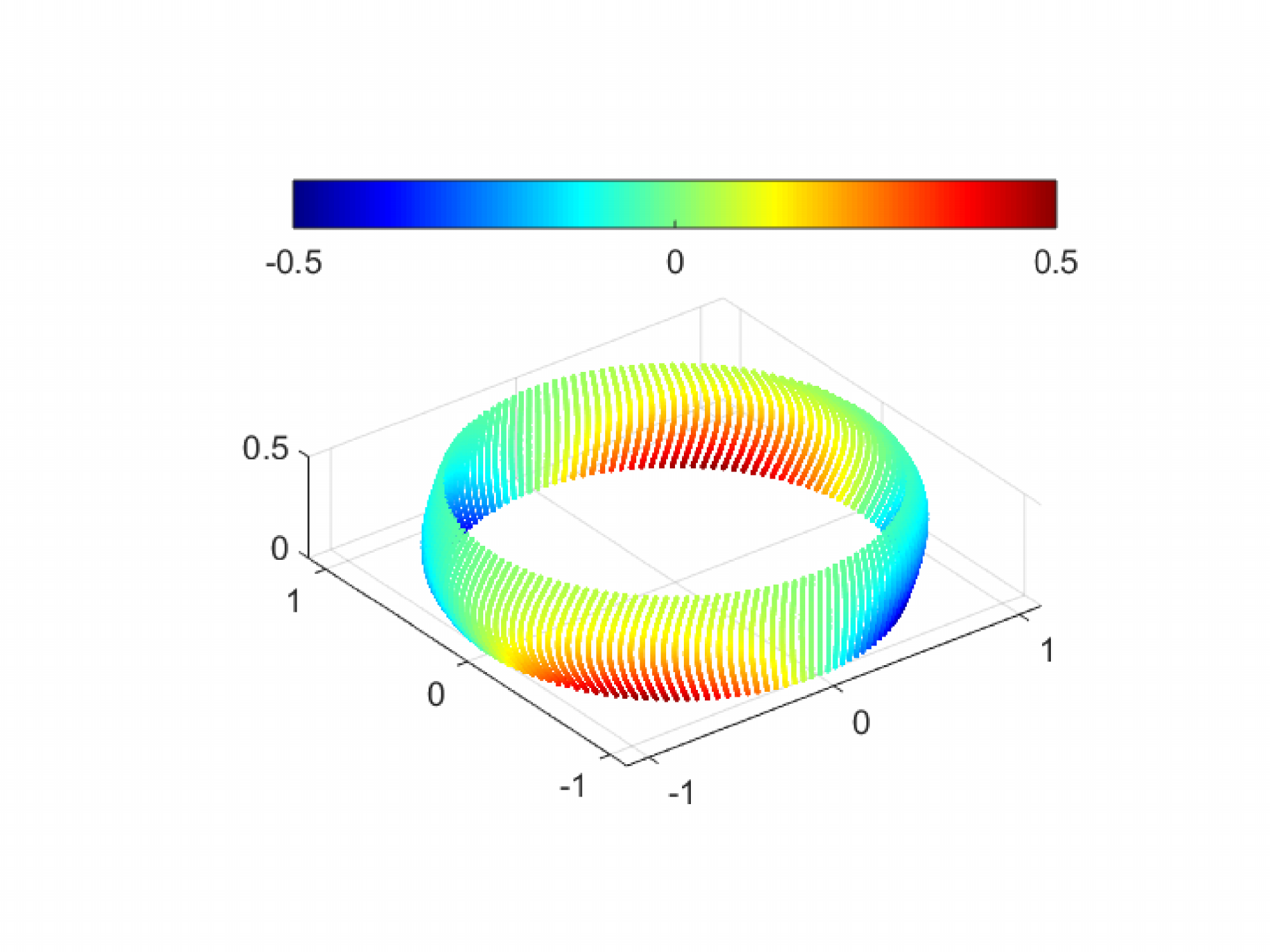}
						\put(2,40) {\scriptsize  $t$}
						\put(80,10) {\scriptsize   $x_1$}
						\put(15,15) {\scriptsize   $x_2$}
						\put(30,0) {\scriptsize  (a) numerical solutions}
					\end{overpic}
					&
					\begin{overpic}
						[width=0.46\textwidth,trim=85 40 85 25, clip=true,tics=10]{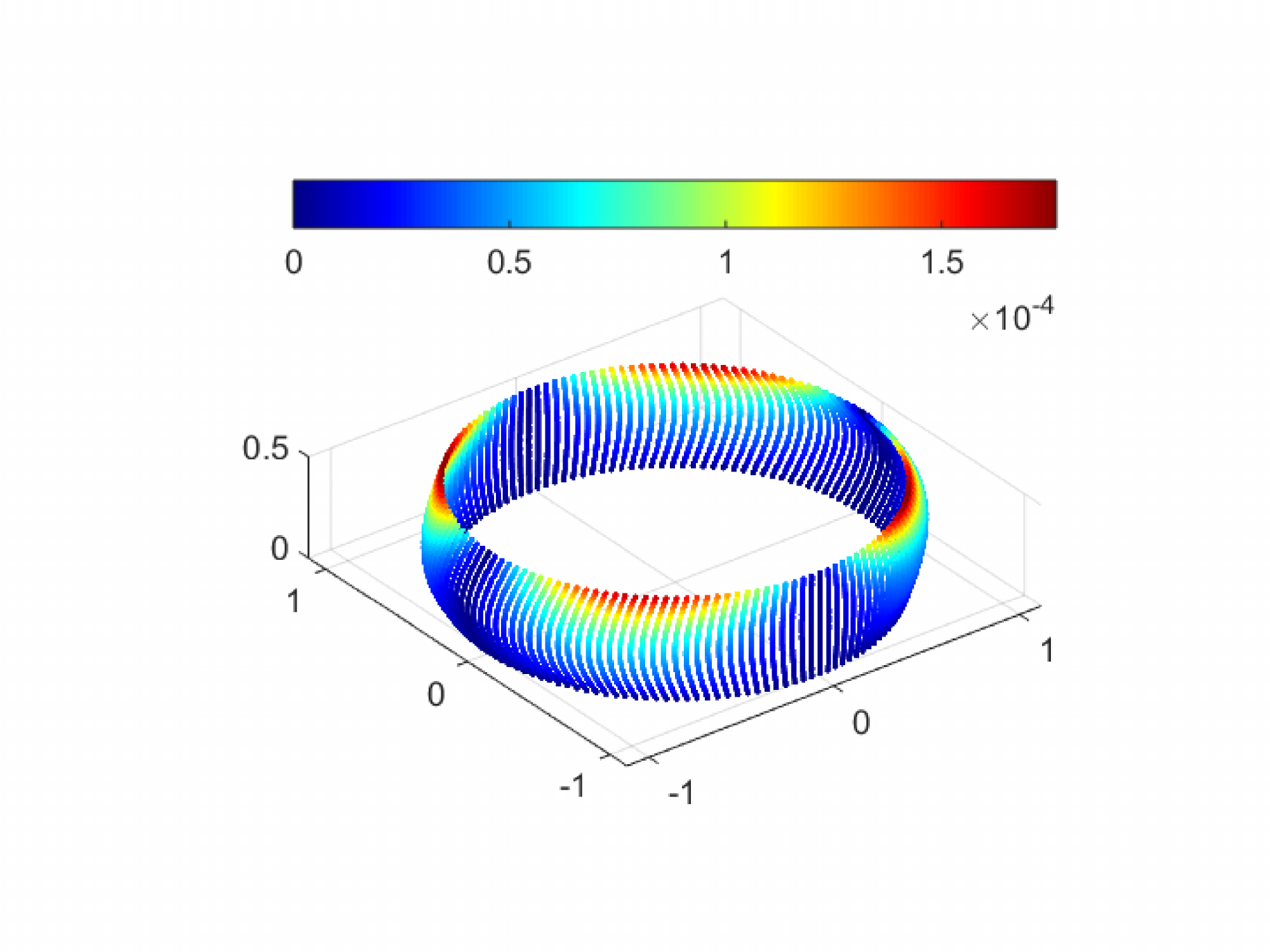}
						\put(30,0) {\scriptsize  (b) error functions}
					\end{overpic}
				\end{tabular}
				\caption{Example~\ref{Ex_MovCur2}. The numerical solutions and error functions (color) relative to $u^*(\vex,t)$ for $t\in[0,0.5]$, obtained by oversampled Algorithm~\ref{KPM-CN} with  $h_Z=2^{-5}\sqrt{2}$ with $(n_X,n_Z)=(214,143)$, using the restricted  Whittle-Mat\'{e}rn-Sobolev kernels of smoothness order $4$, for solving a convective-diffusion equation on a moving curve.}\label{Fig_MovCur}
			\end{figure}

			First, we test the  Algorithm~\ref{KPM-CN} with oversampling.
				In Figure~\ref{Fig_MovCur}, the numerical solutions and their error functions are plotted by color on curves at some time points from $t=0$ to $0.5$, when using $h=2^{-5}\sqrt{2}$.
		 Table~\ref{Tab_MovCur2} shows the corresponding $L^\infty\big([0,T];L^2(\M(t))\big)$ errors under several different oversampling settings of  $n_X>n_Z$  with the corresponding time step $\triangle t =h_Z/4$.
		 With the increasing numbers of both centers and collocation points, the  $L^\infty\big([0,T];L^2(\M(t))\big)$ errors decrease. We then compare results in Table~\ref{Tab_MovCur2} with those listed in Table~\ref{Tab_MovCur} obtained by the corresponding exactly determined settings $X=Z$. We find that using $h_X<h_Z$ can improve accuracy slightly.
Under such smoothness order for the kernel,   oversampling or not both yield the expected second-order $L^\infty\big([0,T];L^2(\M(t))\big)$ convergence. Notably, even for  settings  with very few data points, say $n_Z =5$, we still observe accuracy of $10^{-2}$ in $L^\infty\big([0,T];L^2(\M(t))\big)$ norms that outperforms the FEM mehtod.
	
		\begin{table}
	\centering
	\footnotesize
	\caption{Example~\ref{Ex_MovCur2}. The  $L^\infty\big([0,T];L^2(\M(t))\big)$ errors and  the corresponding estimated order of convergence (eoc) by \eref{eoc} with respect to $h_Z$, obtained by Algorithm~\ref{KPM-CN} with overdetermined formulas ($n_X>n_Z$), under the same other settings as in Figure~\ref{Fig_MovCur}.}\label{Tab_MovCur2}
	\begin{tabular}{l||l|l||cc}
		\hline
		\multirow{1}{*}{$h_Z$}&\multirow{1}{*}{$n_Z$}&\multirow{1}{*}{$n_X$}   
		& $L^\infty\big([0,T];L^2(\M(t))\big)$ & eoc \\
		\hline
		$\sqrt{2}$&$5$& $7$& $2.31823$e-$02$
		&  $-$  \\
		$2^{-1}\sqrt{2}$&$9$ &	$14$ & $ 2.29316$e-$03$
		&  $3.33762$ \\
		$2^{-2}\sqrt{2}$&$18$ &	$27$& $5.44931$e-$04$	&$2.07319$
		\\
		$2^{-3}\sqrt{2}$&$26$ &	$54$& $ 1.34549$e-$04$&	$2.01794$
		\\
		
		\hline
	\end{tabular}
\end{table}

			For comparison with the results of the unfitted FEM as found in \cite{deckelnick2014unfitted},
			Table~\ref{Tab_MovCur} also lists the $L^\infty\big([0,T];L^2(\M(t))\big)$ errors with $\triangle t=2h^2$, as well as the estimates in $L^\infty\big([0,T];\Hsp^2(\M(t))\big)$ error of our Algorithm~\ref{KPM-CN} with $X=Z$. It also presents convergence orders as computed by \eref{eoc}. 
			Algorithm~\ref{KPM-CN}  achieves higher-order accuracy in $L^\infty\big([0,T];L^2(\M(t))\big)$ and $L^\infty\big([0,T];\Hsp^2(\M(t))\big)$  norms, using a relatively large $\triangle t=h/4$.
			The magnitudes of our corresponding $L^\infty\big([0,T];L^2(\M(t))\big)$ errors  are lower by a factor of about 100 relative to the unfitted FEM with many more iterations ($\triangle t=2h^2$).
			
			\begin{table}
				\centering
				\footnotesize
				\caption{Example~\ref{Ex_MovCur2}. The $L^\infty\big([0,T];\Hsp^2(\M(t))\big)$ and $L^\infty\big([0,T];L^2(\M(t))\big)$ errors and the corresponding convergence rates (eoc) with respect to $h$, obtained by Algorithm~\ref{KPM-CN} with exactly determined formulas and the unfitted FEM \cite{deckelnick2014unfitted}, under the same other settings as in Figure~\ref{Fig_MovCur}.}\label{Tab_MovCur}
				\begin{tabular}{l||l|cc|cc|cc}
					\hline
					\multirow{4}{*}{$h$}&\multirow{4}{*}{$n$} & \multicolumn{4}{c|}{Algorithm~\ref{KPM-CN} with $X=Z$ } & \multicolumn{2}{c}{Unfitted FEMs \cite{deckelnick2014unfitted}} \\
					&	&\multicolumn{4}{c|}{ ($\triangle t=h/4$)} & \multicolumn{2}{c}{($\triangle t=2h^2$)} \\
					& &$L^\infty\big([0,T]; $  & \multirow{2}{*}{eoc} & $L^\infty  \big([0,T];$ & \multirow{2}{*}{eoc} & $L^\infty\big([0,T]; $  & \multirow{2}{*}{eoc}\\
						& &\ \, $\Hsp^2(\M(t))\big)$  &   &\ \, $L^2(\M(t))\big)$ &   &\ \, $ L^2(\M(t))\big)$  &  \\
					\hline
					$\sqrt{2}$&$5$ &  $1.14643$e-$1$&  $-$ &     $3.85031$e-$2$                &  $-$ &$-$  &$-$ \\
					$2^{-1}\sqrt{2}$&$9$ &  $8.38177$e-$3$ &  $3.77376$ &       $2.60321$e-$3$               &  $3.88661
					$   & $1.15457$e-$1$ &$-$\\
					$2^{-2}\sqrt{2}$&$18$ &  $2.00492$e-$3$ &  $2.06371$ &       $5.47333$e-$4$               &   $2.24980
					$ &$3.25344$e-$2$ &$1.82732$ \\
					$2^{-3}\sqrt{2}$&$26$ &  $5.23047$e-$4$ &  $1.93853$ &        $1.34576$e-$4$               &   $2.02400
					$ & $8.64172$e-$3$& $1.91258$\\
					$2^{-4}\sqrt{2}$&$72$ &   $1.30694$e-$4$ &  $2.00075$ &       $3.36858$e-$5$               &   $1.99820
					$  & $2.13241$e-$3$& $2.01883$\\
					$2^{-5}\sqrt{2}$&$143$ &   $3.31576 $e-$5$ &  $1.97878$ &       $8.44284$e-$6$               &   $1.99634
					$ & $5.42960$e-$4$ & $1.97357$\\
					\hline
				\end{tabular}
			\end{table}

	\begin{example}\label{Ex_MSp3}{Evolving surface in $\R^3$.}\end{example}

\begin{figure}
	\centering
	\begin{tabular}{cc}
		\begin{overpic}[width=0.46\textwidth,trim=140 60 140 80, clip=true,tics=10]{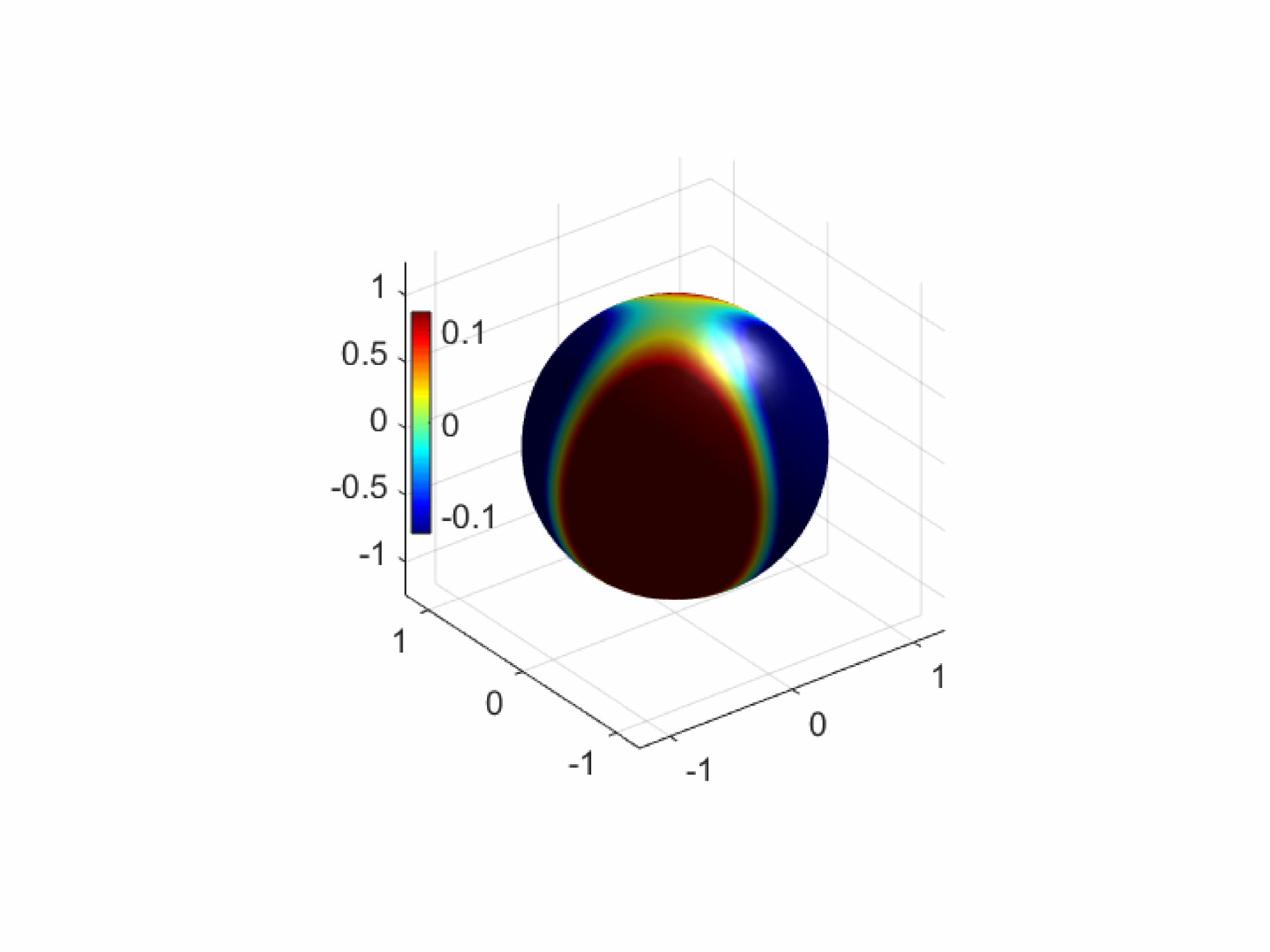}
			\put(47,22) {\scriptsize  (a) $t=0$}
		\end{overpic}
		&
		\begin{overpic}[width=0.46\textwidth,trim=140 60 140 80, clip=true,tics=10]{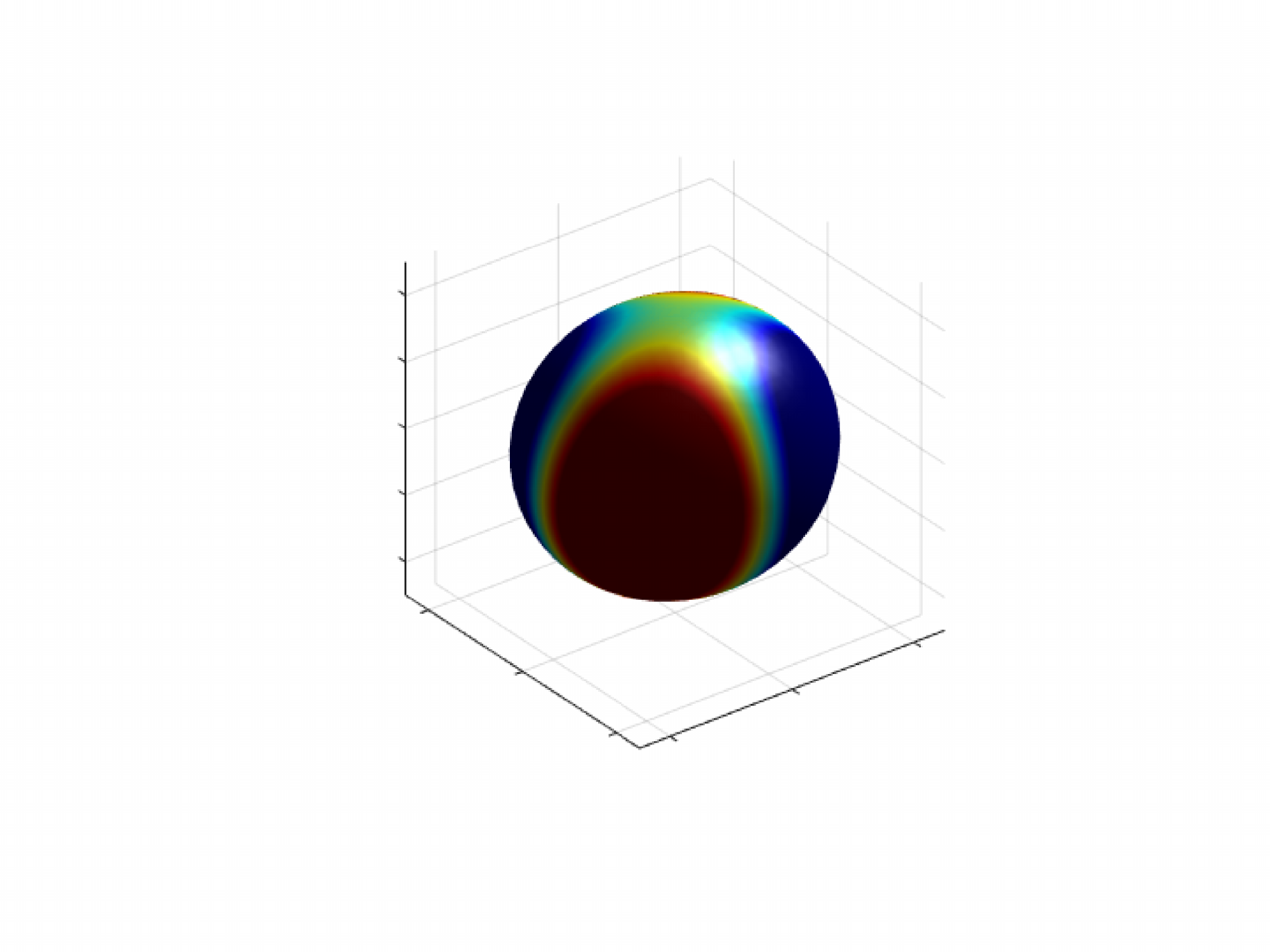}
			\put(47,22) {\scriptsize  (b) $t=0.5$}
		\end{overpic}\\
		\begin{overpic}[width=0.46\textwidth,trim=140 60 140 80, clip=true,tics=10]{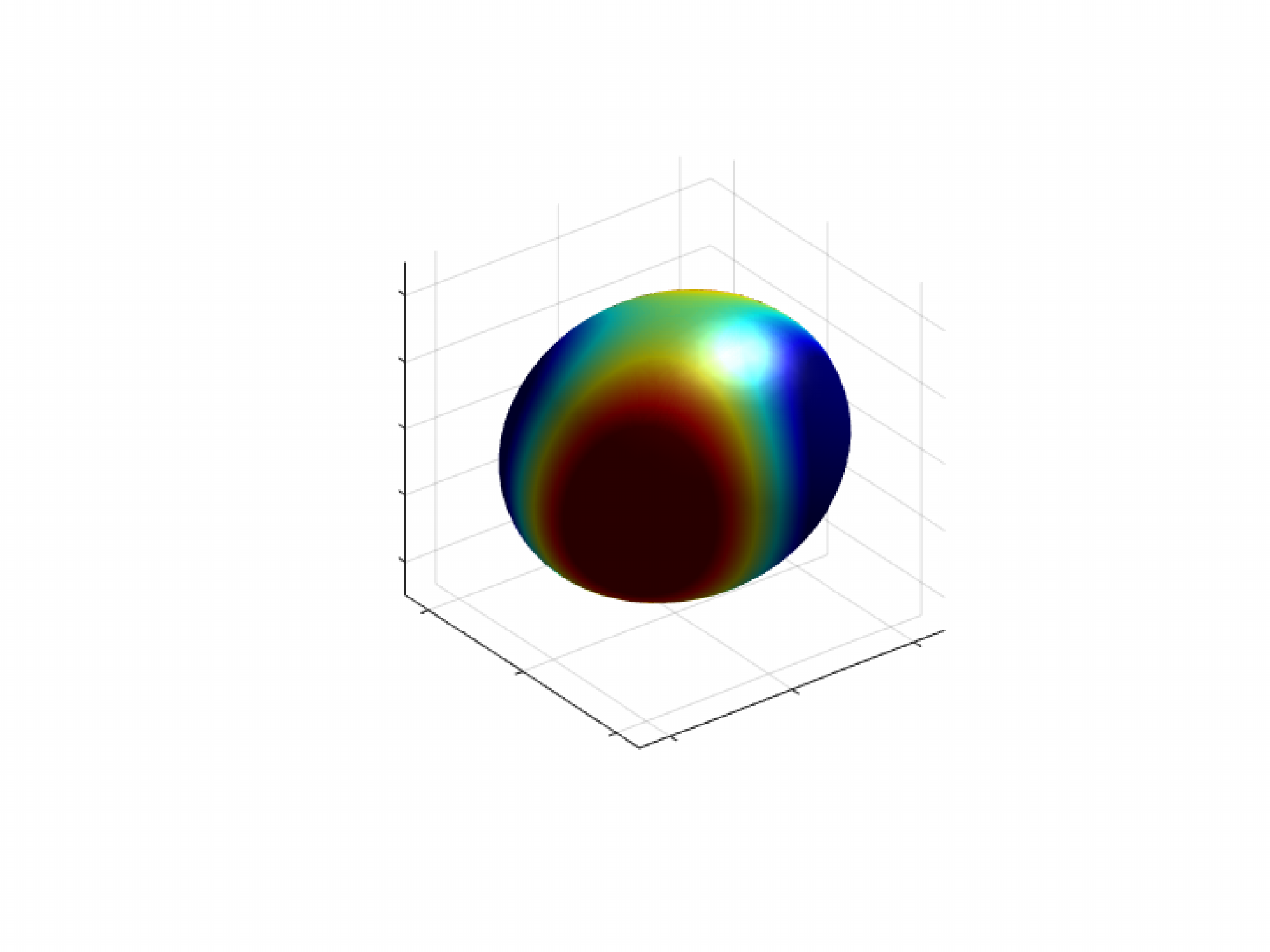}
			\put(47,22) {\scriptsize  (c) $t=1.125$}
		\end{overpic}
		&\begin{overpic}[width=0.46\textwidth,trim=140 60 140 80, clip=true,tics=10]{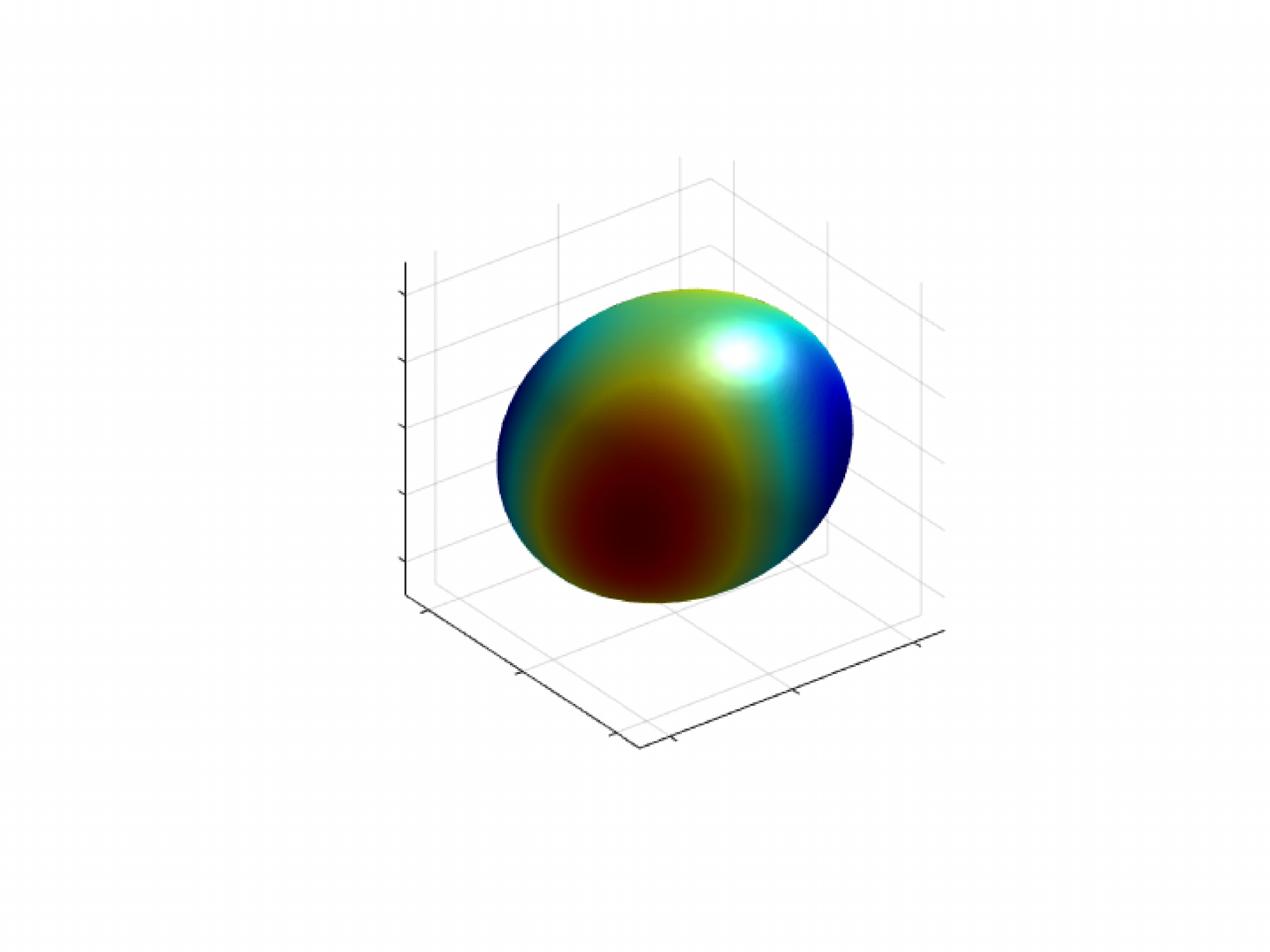}
			\put(47,22) {\scriptsize  (d) $t=1.75$}
		\end{overpic}
		\\
		\begin{overpic}[width=0.46\textwidth,trim=140 90 140 80, clip=true,tics=10]{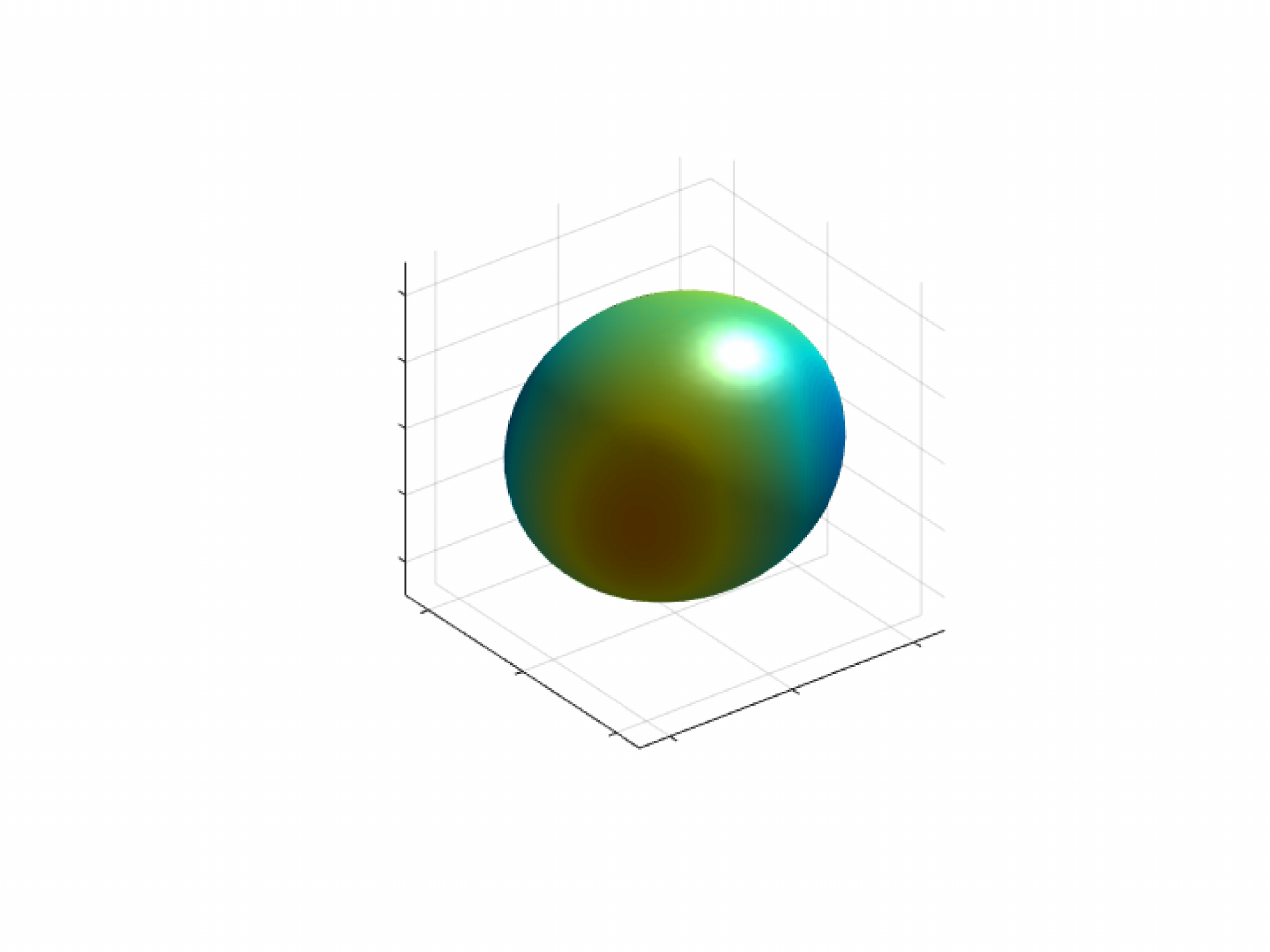}
			\put(47,12) {\scriptsize  (e) $t=2.375$}
		\end{overpic}
		&
		\begin{overpic}[width=0.46\textwidth,trim=140 90 140 80, clip=true,tics=10]{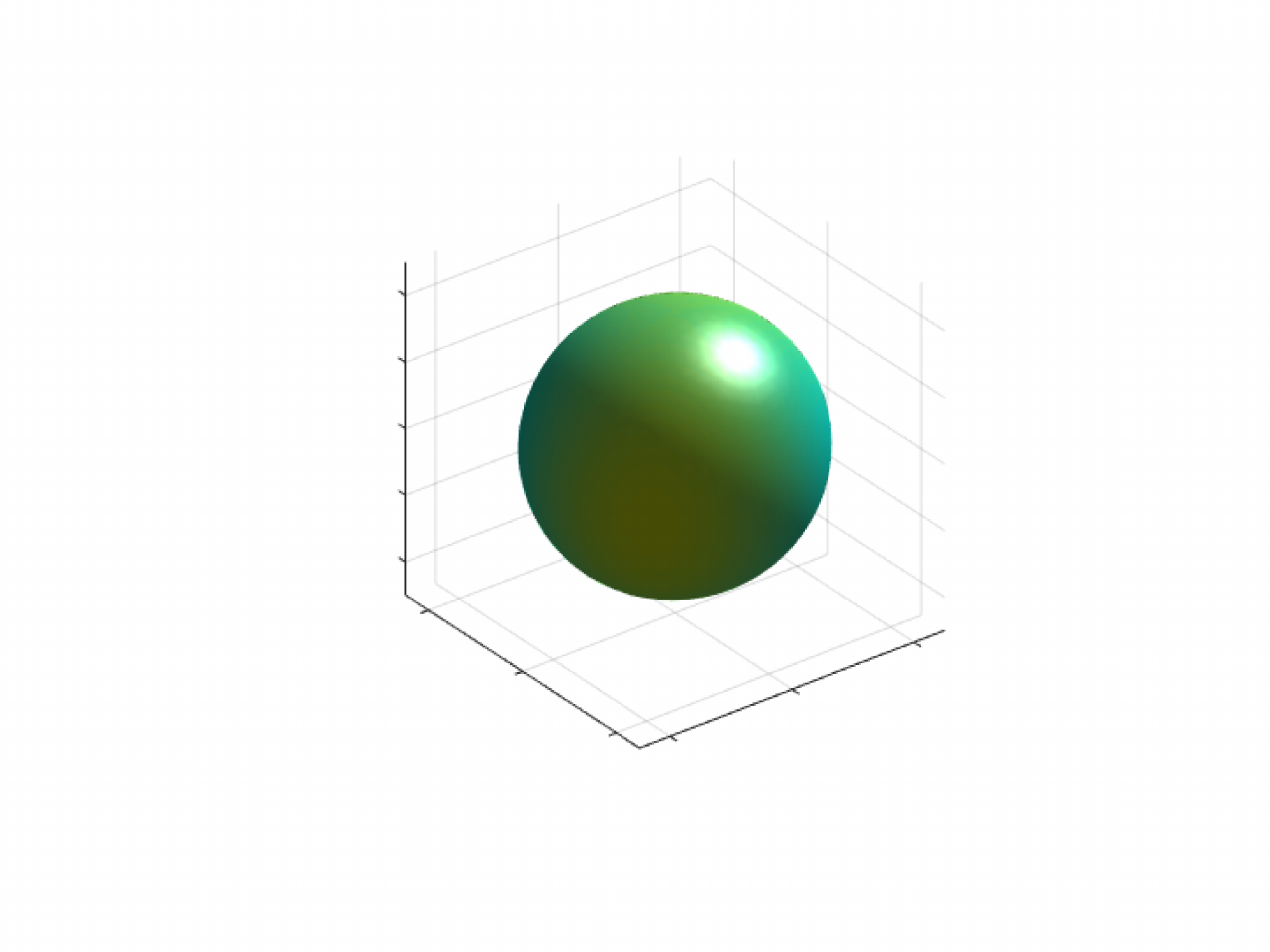}
			\put(47,12) {\scriptsize  (f) $t=3$}
		\end{overpic}
	\end{tabular}
	\caption{Example~\ref{Ex_MSp3}. Snapshots of numerical solutions (color) at several times, obtained by  Algorithm~\ref{KPM-CN}, using the kernel of smoothness order $4$ and $h=\triangle t =1/8$,
on a moving surface (see a movie from { https://youtu.be/GPo-mZaDCuE}).}\label{Fig_MovSph1}
\end{figure}

This is a 3D version of the previous example, but we focus on the exactly determined case with $X^m=Z^m$ at all time steps.
We use kernels with the smoothness order $4$.
		 In \cite{grande2014eulerian}, two Eulerian FEMs (the cGdG and  cGcG methods) employed Galerkin methods both in time and space with narrow bands to solve \eref{eq_CDpupt}.
		The test problem in \cite[Section~11.3]{grande2014eulerian} is posed on a moving surface 
			given by
			\[
			\M(\vex,t)=\Big(\frac{x_1}{1+0.25\sin(t)}\Big)^2+x_2^2+x_3^2-1=0\ \mbox{\quad for } t\in[0,4].
			\]
We solve the same problem by Algorithm~\ref{KPM-CN} with $\triangle t =1/8$.
Based on the initial points with $h:=h_{Z^0}=\triangle t$, we move each point by point-to-point motion based on the analytical formula in polar coordinate, that is, the fill distance $h_{Z^m}$ slightly varies in time periodically indeed.
In Figure~\ref{Fig_MovSph1}, we present a few snapshots of the evolved surfaces and numerical solutions by colormap.
		
\begin{figure}
				\centering
				\begin{tabular}{cc}
					\begin{overpic}[width=0.46\textwidth,trim=85 40 85 25, clip=true,tics=10]{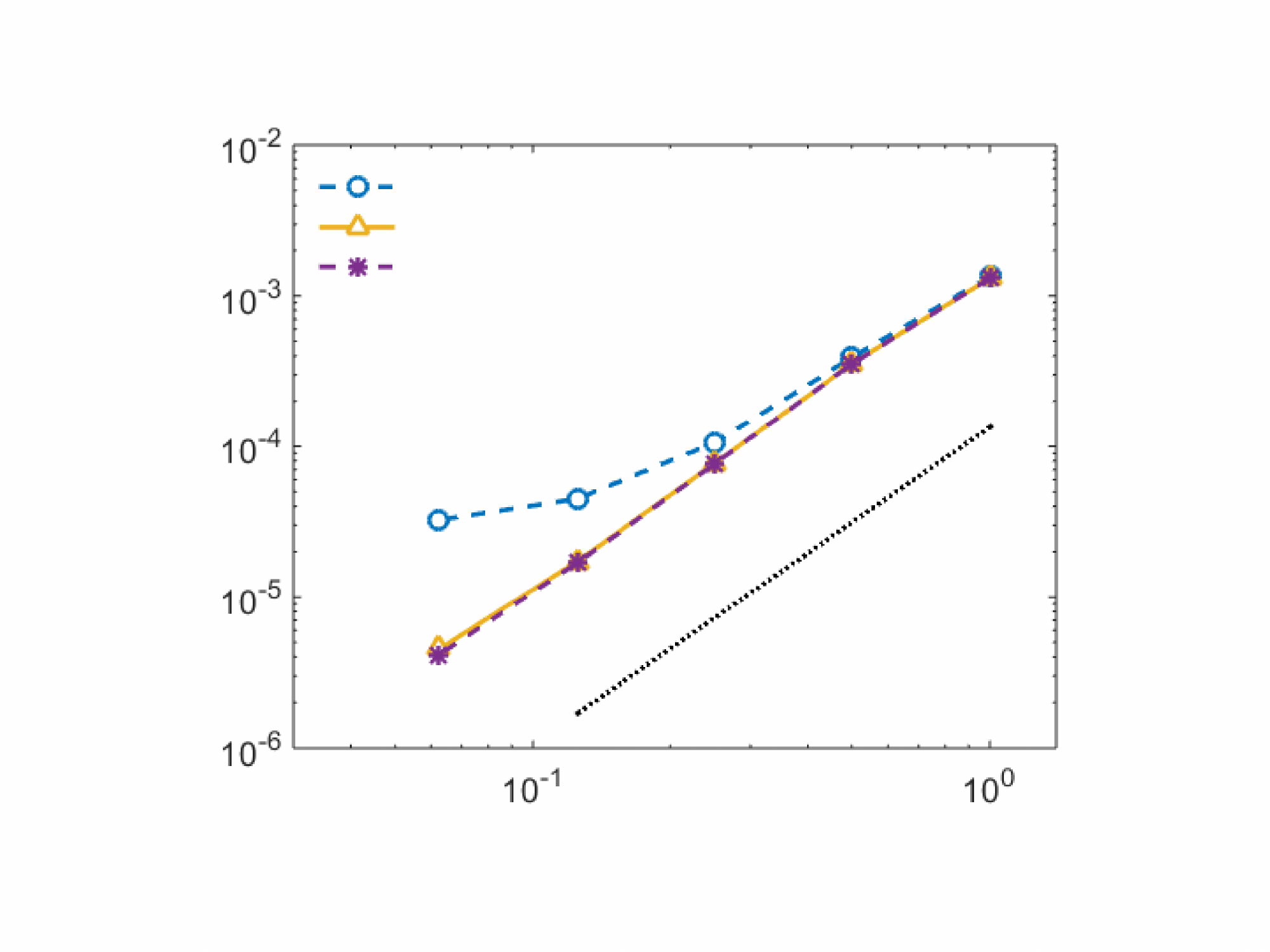}
						\put(-4,20) {\scriptsize \rotatebox{90}{$L^\infty\big([0,T];L^2(\M(t))\big)$ errors}}
						\put(95,7) {\scriptsize   $\triangle t$}
						
						\put(25,-2) {\scriptsize  (a) convergence for $\triangle t$}
						
						\put(23.5,74) {\scriptsize  $h=1/2$}
						\put(23.5,69.5) {\scriptsize  $1/4$}
						\put(23.5,65) {\scriptsize  $1/8$}
						\put(72,22) {\scriptsize  $2.1$}
					\end{overpic}
					&
					\begin{overpic}[width=0.46\textwidth,trim=85 40 85 25, clip=true,tics=10]{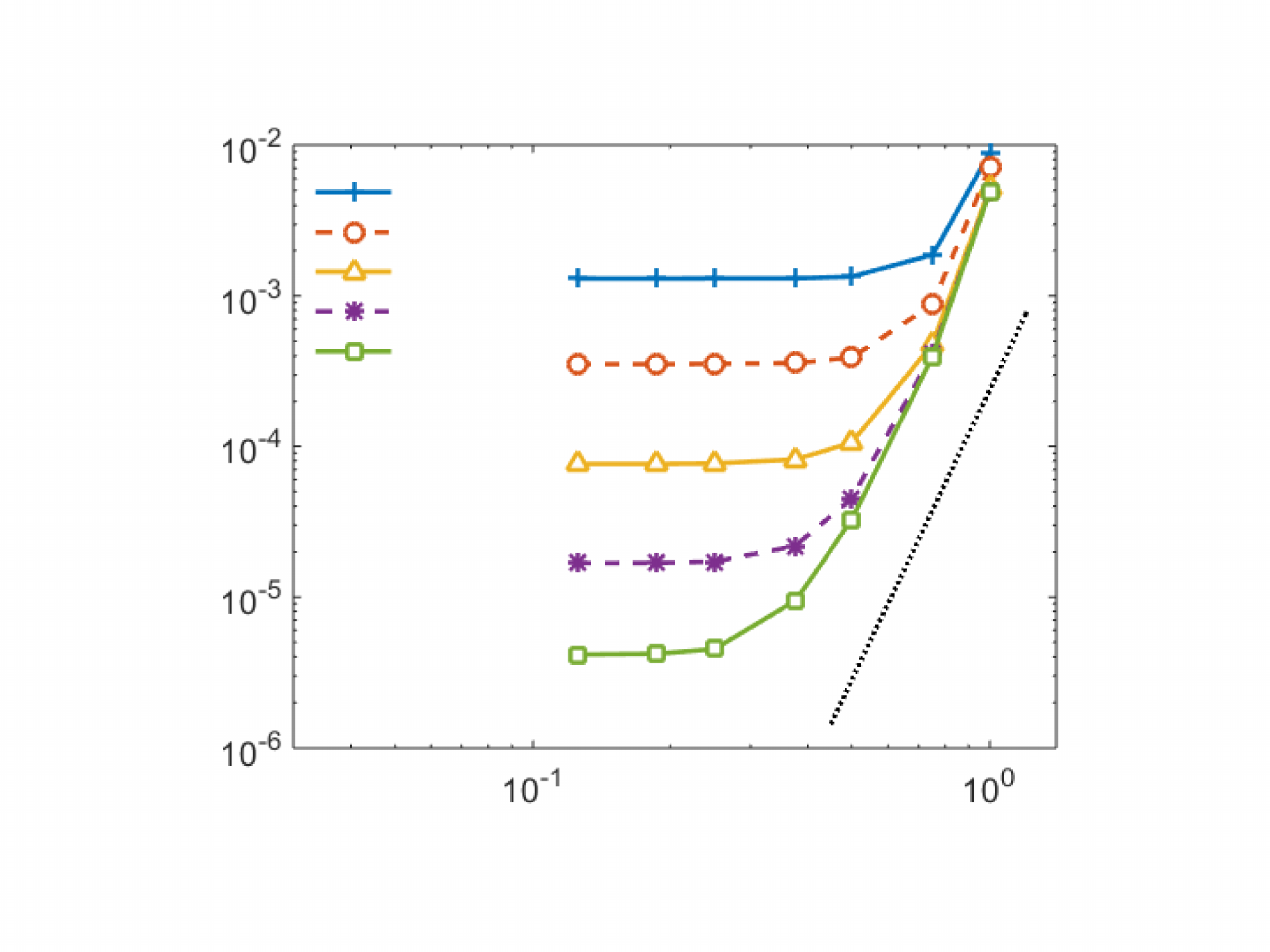}
						\put(95,7) {\scriptsize   $h$}
						\put(25,-2) {\scriptsize  (b) convergence for $h$}
						
						\put(23.5,74) {\scriptsize  $\triangle t=1$}
						\put(23.5,69.5) {\scriptsize  $1/2$}
						\put(23.5,65) {\scriptsize  $1/4$}
						\put(23.5,60.5) {\scriptsize  $1/8$}
						\put(23.5,56) {\scriptsize  $1/16$}
						\put(80,22) {\scriptsize  $6.4$}
					\end{overpic}
				\end{tabular}
				\caption{Example~\ref{Ex_MSp3}. The $L^\infty\big([0,T];L^2(\M(t))\big)$ error profiles with different $h$ and $\triangle t$, obtained by our proposed Algorithm~\ref{KPM-CN}, under the same settings as in Figure~\ref{Fig_MovSph1}.}\label{Fig_MovSph}
			\end{figure}

			To be compared with \cite[Figures~13-14]{grande2014eulerian}, Figure~\ref{Fig_MovSph}  illustrates the convergence profile of Algorithm~\ref{KPM-CN} in   $L^\infty\big([0,T];L^2(\M(t))\big)$ errors.		
First, Figure~\ref{Fig_MovSph}(a) is obtained by three fixed $h$-values as $\triangle t$ varies.  It can be seen that errors for $h=1/4$ are almost as accurate as those for $h=1/8$, under the same $\triangle t$.

Our comparison with the cGdG and cGcG approaches in \cite{grande2014eulerian} is as follows.
With respect to time, our convergence rates and those of the two FEMs are all around order 2 due to CN. In terms of accuracy, we  obtain errors  less than $10^{-5}$ by using $h=1/4$, while their best FEM result shown is approximately $10^{-4}$ using a much smaller $h=1/64$.
For spatial convergence, we fix a several $\triangle t$  and let $h$ is allowed to vary to obtain Figure~\ref{Fig_MovSph}(b).
Generally speaking, our numerical solutions are more accurate based on \cite[Figures~13-14]{grande2014eulerian}, which used smaller $h$-values for each $\triangle t$.
Our method shows spatial convergence of $6.4$, whereas that  in \cite{grande2014eulerian} is around $2$.

			Table~\ref{Tab_MovSph} further lists the $L^\infty\big([0,T];\Hsp^2(\M(t))\big)$ and $L^\infty\big([0,T];\Hsp^1(\M(t))\big)$ estimates obtained by Algorithm~\ref{KPM-CN}, and the $L^\infty\big([0,T];\Hsp^1(\M(t))\big)$ results of two Eulerian FEMs in \cite{grande2014eulerian} using various $h=\triangle t$.
Note that  Algorithm~\ref{KPM-CN} remains second order while the FEM methods drop to first order.
			\begin{table}
				\centering
			\scriptsize
				\caption{Example~\ref{Ex_MSp3}. The corresponding $L^\infty\big([0,T];\Hsp^2(\M(t))\big)$ and $L^\infty\big([0,T];\Hsp^1(\M(t))\big)$ errors and  the estimated order of convergence  (eoc) by \eref{eoc} using $h=\triangle t$, compared with the $L^\infty\big([0,T];\Hsp^1(\M(t))\big)$ results in \cite[Table~4]{grande2014eulerian}, under the same settings as in Figure~\ref{Fig_MovSph1}.}\label{Tab_MovSph}
				\begin{tabular}{l||l|cc|cc|cc|cc}
					\hline
					\multirow{3}{*}{$h$}&\multirow{3}{*}{$n $} & \multicolumn{4}{c|}{Algorithm~\ref{KPM-CN}} & \multicolumn{2}{c}{cGdG} &\multicolumn{2}{c}{cGcG} \\ 
					& & $L^\infty\big([0,T];$ & \multirow{2}{*}{eoc} & $L^\infty\big([0,T];$ & \multirow{2}{*}{eoc} & $L^\infty\big([0,T];$ & \multirow{2}{*}{eoc} & $L^\infty\big([0,T];$ & \multirow{2}{*}{eoc}\\
					& &\ \ $\Hsp^2(\M(t))\big)$ &   & \ \ $\Hsp^1(\M(t))\big)$ &   & \ \ $\Hsp^1(\M(t))\big)$ &   &\ \ $\Hsp^1(\M(t))\big)$ &  \\\hline
					$3/4$&$36$& $5.5926$e-$3$ & $-$ & $3.7442$e-$3$               &  $-$       & $-$ & $-$ & $-$&$-$\\
					$1/2$&$78$&    $1.9828$e-$3$ & $2.557$ & $9.2059$e-$4$               &  $3.460$       & $0.574872$ & $-$ & $1.12954$&$-$\\
					$1/4$&$312$&  $4.3551$e-$4$ & $2.187$ &  $1.8276$e-$4$               &   $2.333$  & $0.303749$ & $0.920$&$0.550169$&$1.04$\\
					$1/8$&$1206$&     $9.6060$e-$5$ & $2.181$ &  $4.0135$e-$5$               &   $2.187$  & $0.157731$ & $0.945$&$0.249556$& $1.14$\\
					$1/16$&$4836$&    $2.3407$e-$5$ & $2.037$ &  $9.7715$e-$6$               &   $2.038$  & $0.0803239$ & $0.974$ &$0.112138$&$1.15$\\
					\hline
				\end{tabular}
			\end{table}



	\begin{example}\label{Ex_masscons}{Mass conservation on an expanding sphere.}\end{example}
		
Aiming to verify the mass conservation law in \eref{eq_conlaw1}, we solve PDE \eref{eq_CDpupt} without the source term, i.e., $f=0$.
We choose a small diffusion coefficient $\df=10^{-3}$ to approximate the viscosity solutions to \eref{eq_CDpupt}. The tested surface is defined by
			\[
			\M(\vex,t)=x_1^2+x_2^2+x_3^2-(e^{t/5})^2=0\ \mbox{\quad for } t\in[0,1].
			\]
		 We set the initial solution to be $u^*(\vex,0)=0.5+x_1 x_2 x_3$.
			
\begin{figure}
				\centering
				\begin{overpic}
					[width=0.32\textwidth,trim=190 120 180 100, clip=true,tics=10]{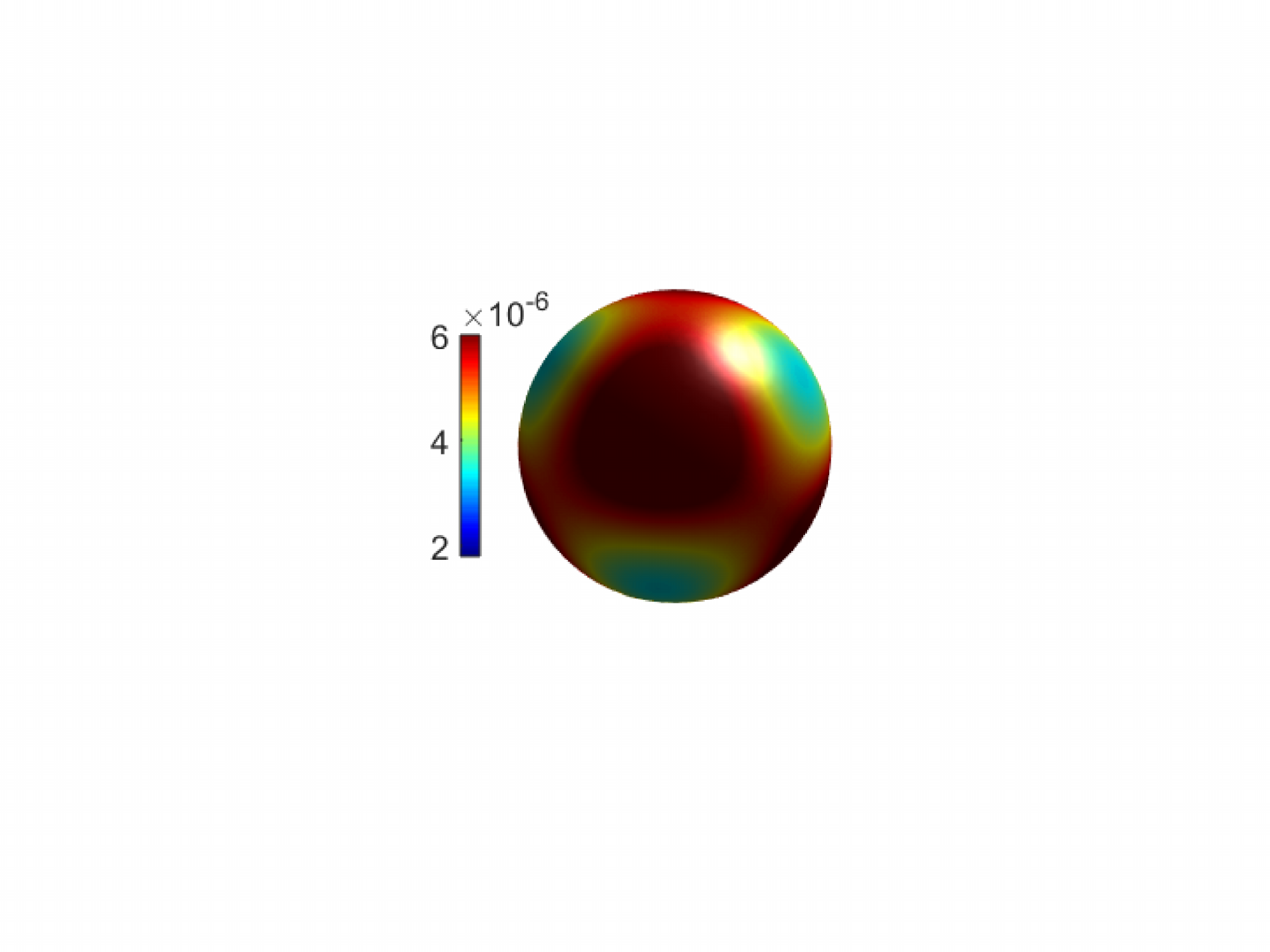}
					\put(42,0) {\scriptsize  (a) $t=0.3$}
				\end{overpic}
				\begin{overpic}
					[width=0.32\textwidth,trim=190 120 180 100, clip=true,tics=10]{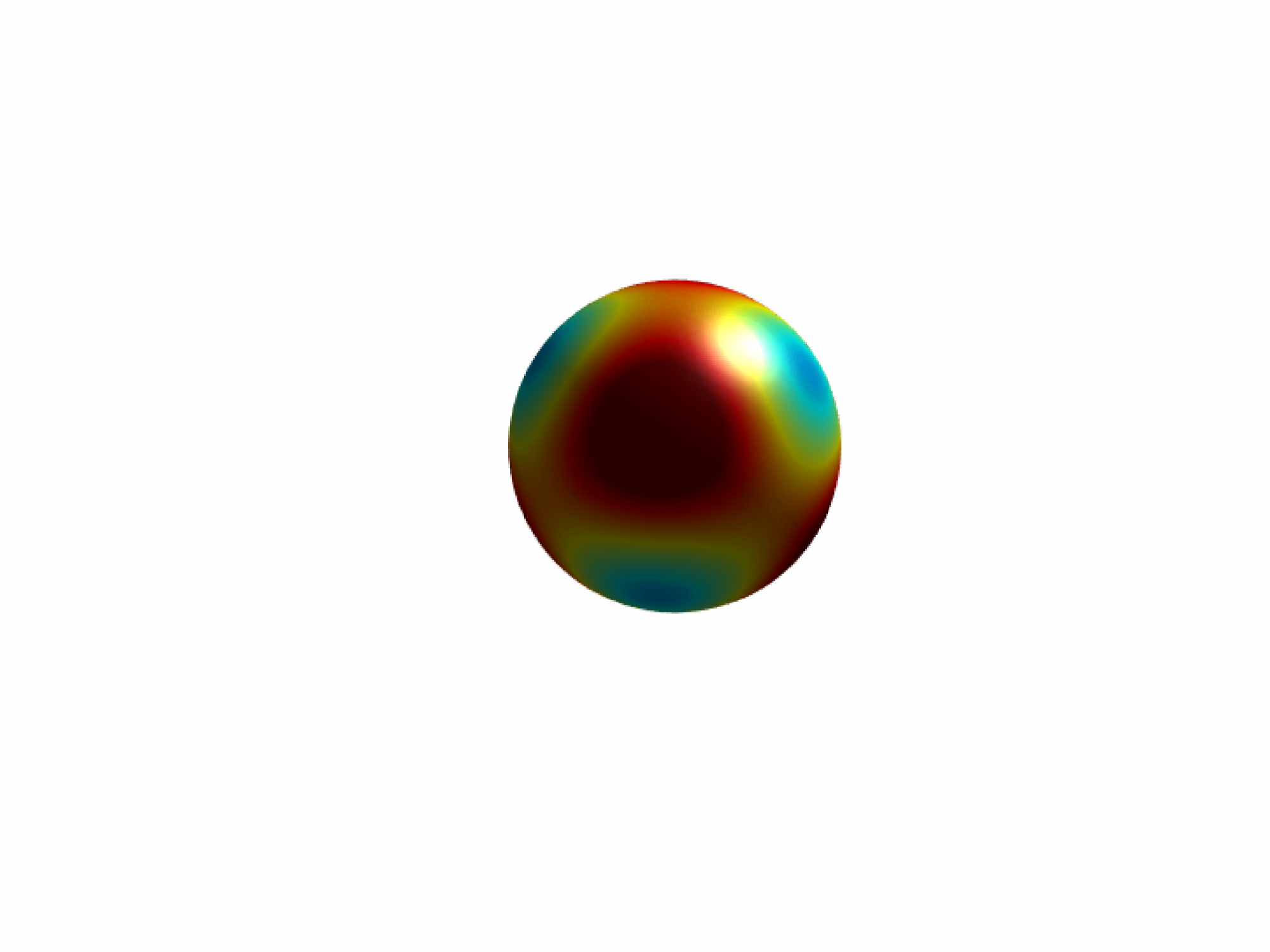}
					\put(49,0) {\scriptsize  (b) $t=0.6$}
				\end{overpic}
				\begin{overpic}
					[width=0.32\textwidth,trim=190 120 180 100, clip=true,tics=10]{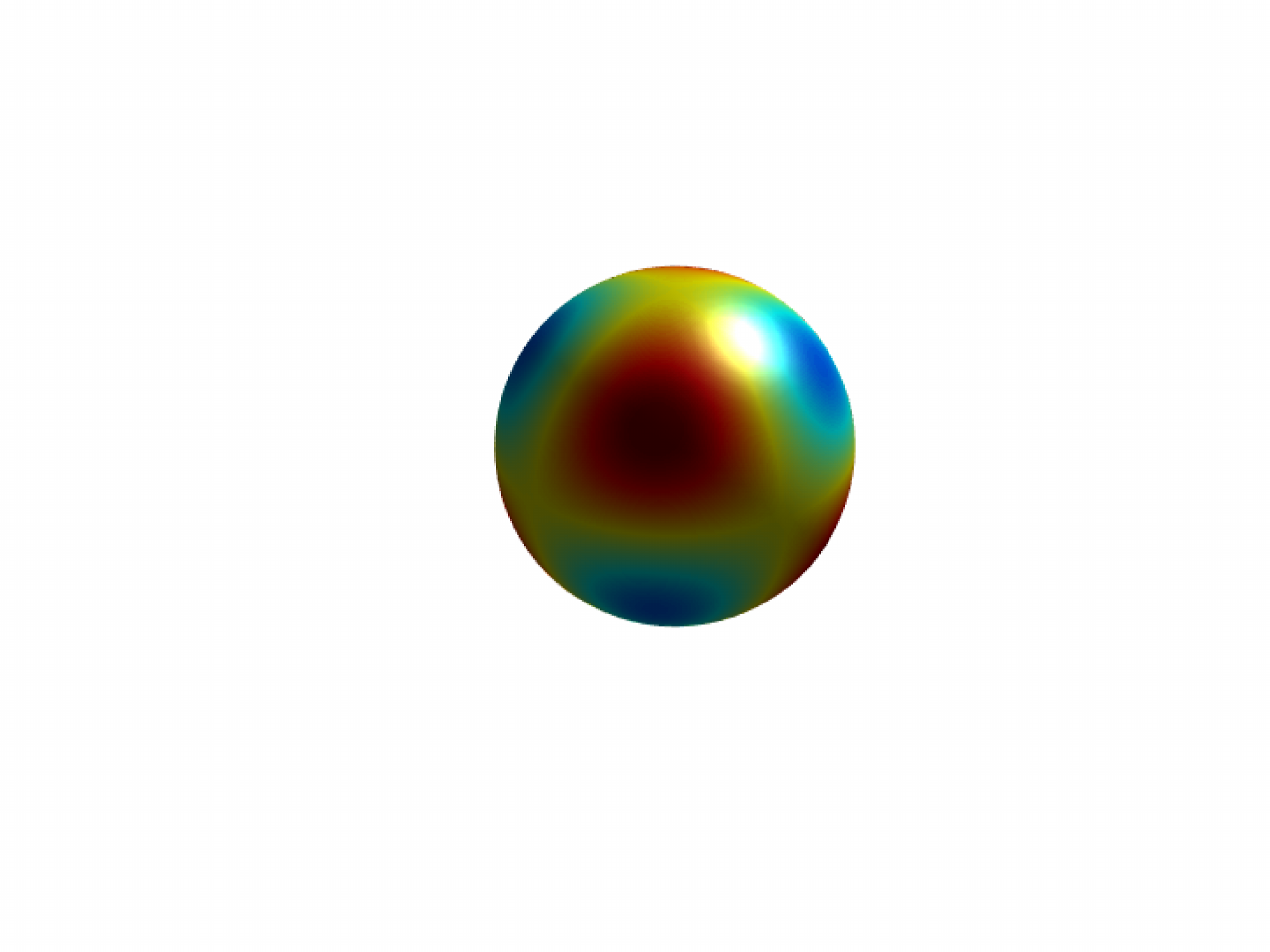}
					\put(49,0) {\scriptsize  (c) $t=1$}
				\end{overpic}
				\caption{Example~\ref{Ex_masscons}. Difference functions of solutions obtained by Algorithms~\ref{KPM-CN} and \ref{AKM-CN}a, under $h=\triangle t =0.02$, for solving a convective-diffusion equation on a sphere expanding along outer normal directions.}\label{Fig_SphMass1_er}
			\end{figure}

We  use $X=Z$ in both Algorithms~\ref{KPM-CN} and \ref{AKM-CN}
(Algorithm~\ref{AKM-CN}a and b are practically identical  here)
to solve this problem with the kernels of order $4$ and $h=\triangle t =0.02$.
Solutions of two algorithm agree up to some order of $10^{-6}$;
their difference on $\M(t)$ is shown in Figure~\ref{Fig_SphMass1_er} for a particular parameter setup.

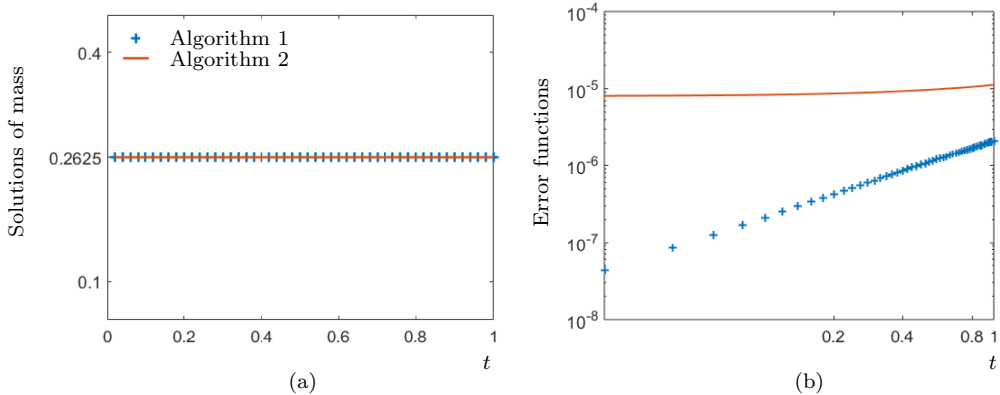
\begin{figure}
				\centering
				\begin{tabular}{cc}
					\begin{overpic}[width=0.46\textwidth,trim=70 17 85 30, clip=true,tics=10]{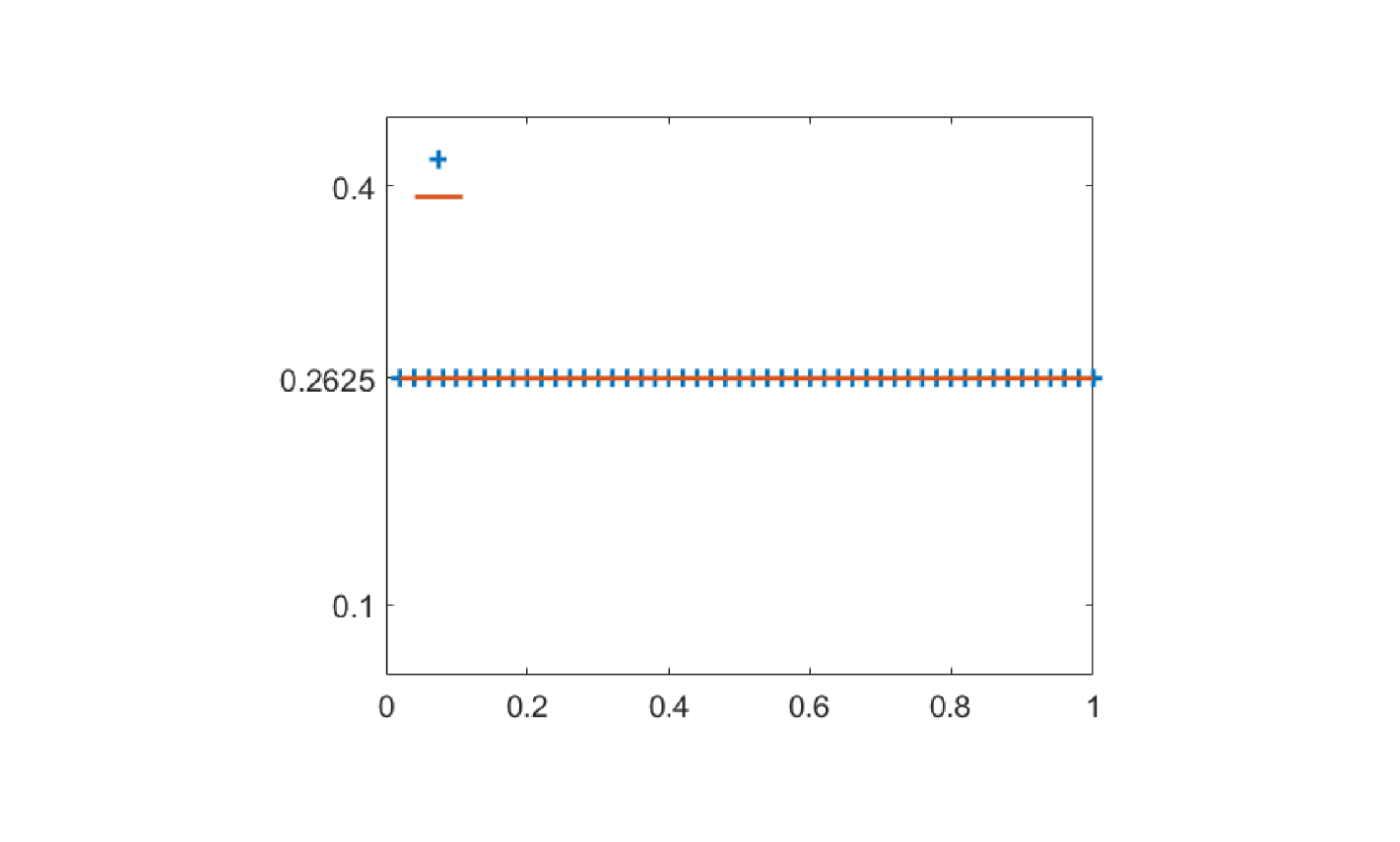}
						
						\put(-4,30) {\scriptsize \rotatebox{90}{Solutions of mass}}
						
						\put(96,2) {\scriptsize   $t$}
						\put(55,-2) {\scriptsize   (a)}
						\put(30,70.5) {\scriptsize   Algorithm~\ref{KPM-CN}}
						\put(30,66) {\scriptsize   Algorithm~\ref{AKM-CN}}
					\end{overpic}
					&	\begin{overpic}[width=0.46\textwidth,trim=70 40 85 25, clip=true,tics=10]{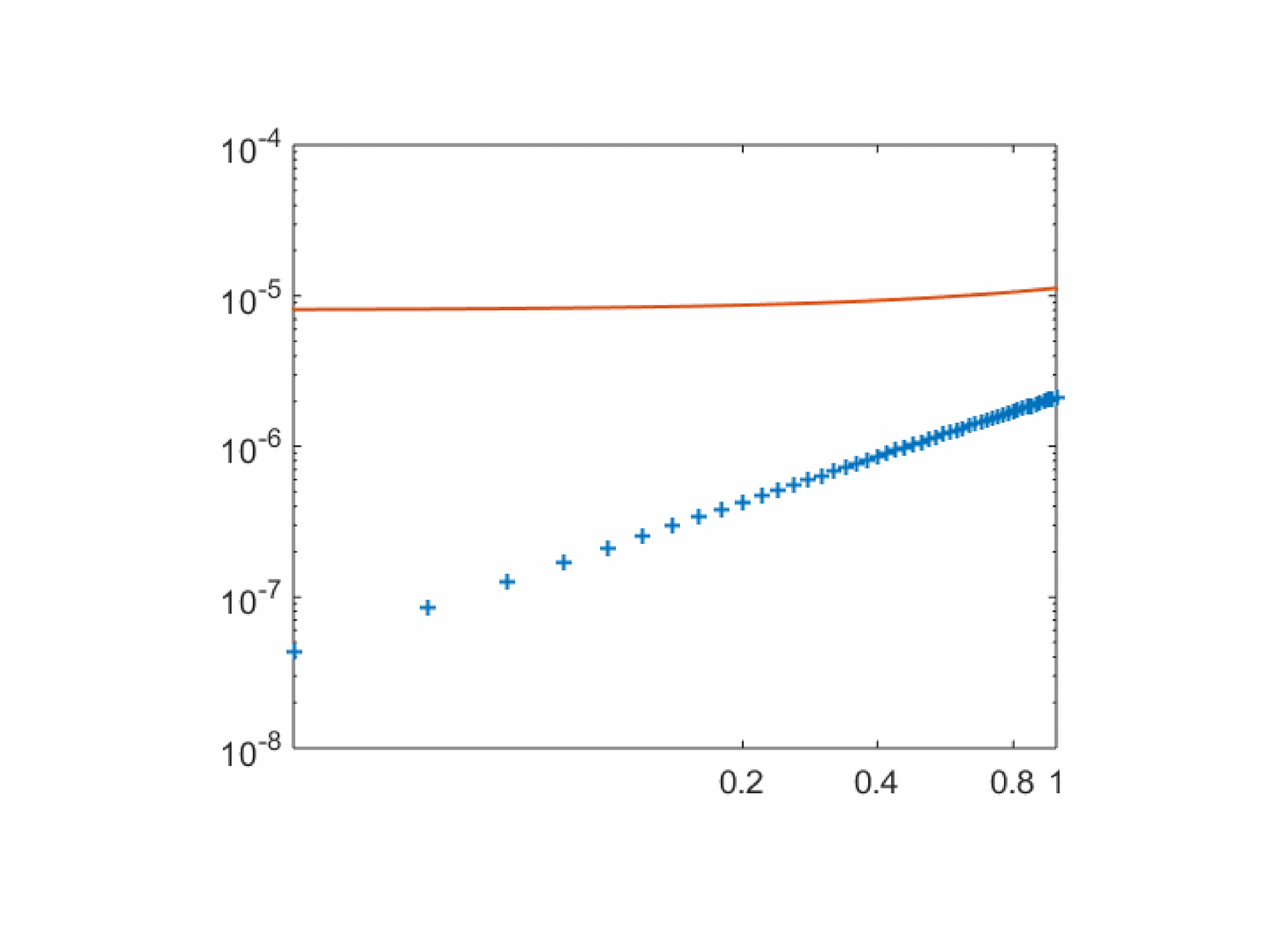}
						
						\put(-0,32) {\scriptsize \rotatebox{90}{Error functions}}
						\put(96,2) {\scriptsize   $t$}
						\put(55,-2) {\scriptsize   (b)}
						
					\end{overpic}
				\end{tabular}
				\caption{Example~\ref{Ex_masscons}. (a) Solutions of mass and (b) error functions between exact mass and  them, obtained by Algorithms~\ref{KPM-CN} and \ref{AKM-CN} respectively, under the same settings as in Figure~\ref{Fig_SphMass1_er}.}\label{Fig_SphMass2}
			\end{figure}

Now we can see Figure~\ref{Fig_SphMass2} for the masses of the solutions
$\int_{\M(t)}u\ \text{d}\M$ over time,
which are approximated by
$\sum_{i}^{n_X} u(\vex_i)\triangle \M_i$ at $X$ based on some triangle mesh $\mathcal{T}$. It is easy to see that $u(\vex_i)\triangle \M_i=u(\vex_i)\mbox{Area}(\mathcal{T}_i)/3$ where $\mathcal{T}_i$ is the triangle containing $\vex_i$.
For Figure~\ref{Fig_SphMass2}(a), it is clear that both algorithms conserve mass as desired. When we subtract the exact mass computed using $u^*(\vex,0)$ to yield errors in Figure~\ref{Fig_SphMass2}(b), we can see the different mass-conservation behaviour of the proposed algorithms. The analytically exact Algorithm~\ref{KPM-CN} shows gradually increasing error in time, which is typical for accumulative error. On the other hand, the nearly-constant error curve of Algorithm~\ref{AKM-CN}  in Figure~\ref{Fig_SphMass2}(b) suggests that the error in numerical mass comes from the approximation steps in formulation.

			
%

\begin{example}\label{Ex_ContSurf}{Continuous initial conditions on point cloud with high curvatures.}\end{example}

The last two examples aim to prepare for problems on merging surfaces.
In both cases, we consider the mean-curvature motion of a kissing-spheres, as in Figure \ref{Fig_MovMeanCur}, which has high initial curvature at the contact point. Here, we use the FDM proposed in \cite{bronsard1995numerical} to obtain point clouds in $\R^3$ at each discrete time.
Without analytic formula for $\M(t)$, we can only use Algorithm~\ref{AKM-CN}.

%

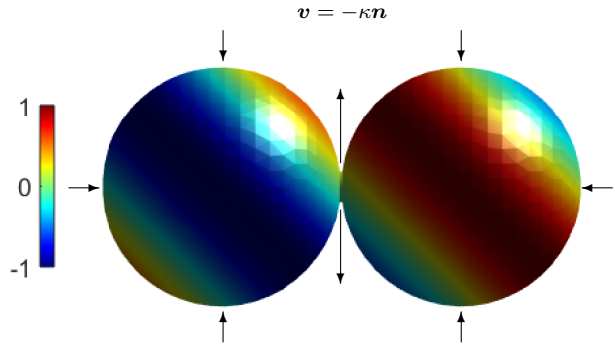
\begin{figure}
	\centering	
	\begin{overpic}[width=0.6\textwidth,trim=90 80 90 50, clip=true,tics=10]{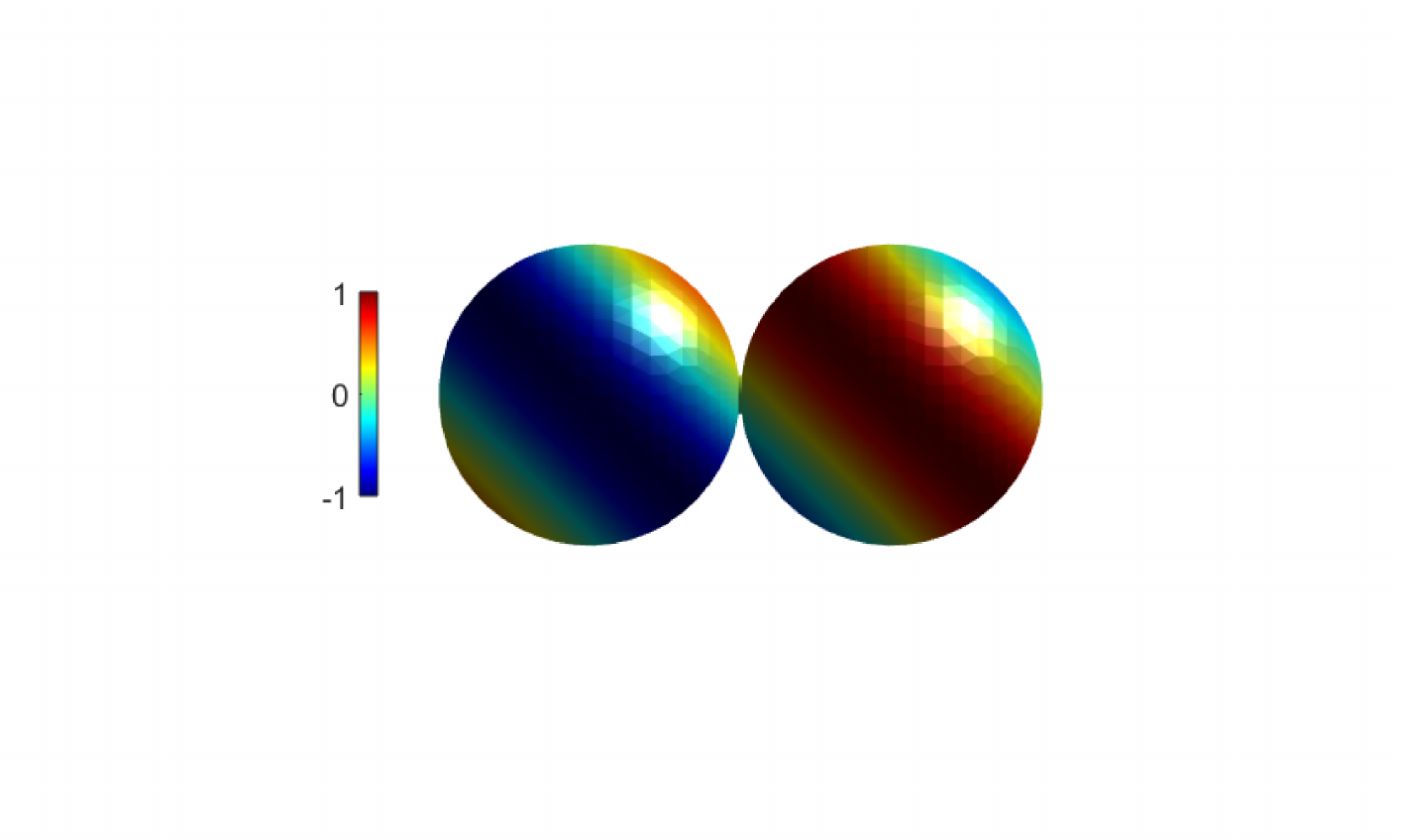}
		\put(55,25.5){\vector(0,1){12}}
		\put(55,18){\vector(0,-1){12}}
		
		\put(36,47){\vector(0,-1){5}}
		\put(74.5,47){\vector(0,-1){5}}
		\put(36,-3.5){\vector(0,1){5}}
		\put(74.5,-3.5){\vector(0,1){5}}
		
		\put(11,21.5){\vector(1,0){5}}
		\put(99,21.5){\vector(-1,0){5}}
		
		\put(48,49){\scriptsize $\vev=-\kappa \ven$}
	\end{overpic}
	\caption{Example~\ref{Ex_ContSurf}. A continuous IC (color) on a surface that evolves by mean curvature motion (arrows are the velocities $\vev$).}\label{Fig_Contin}
\end{figure}

\begin{figure}
	\centering	
	\begin{tabular}{cc}
		\multicolumn{2}{c}{\scriptsize Exactly determined formula with $X=Z$ ($n=1196$)}\\
		\begin{overpic}[width=0.46\textwidth,trim=100 70 102 60, clip=true,tics=10]{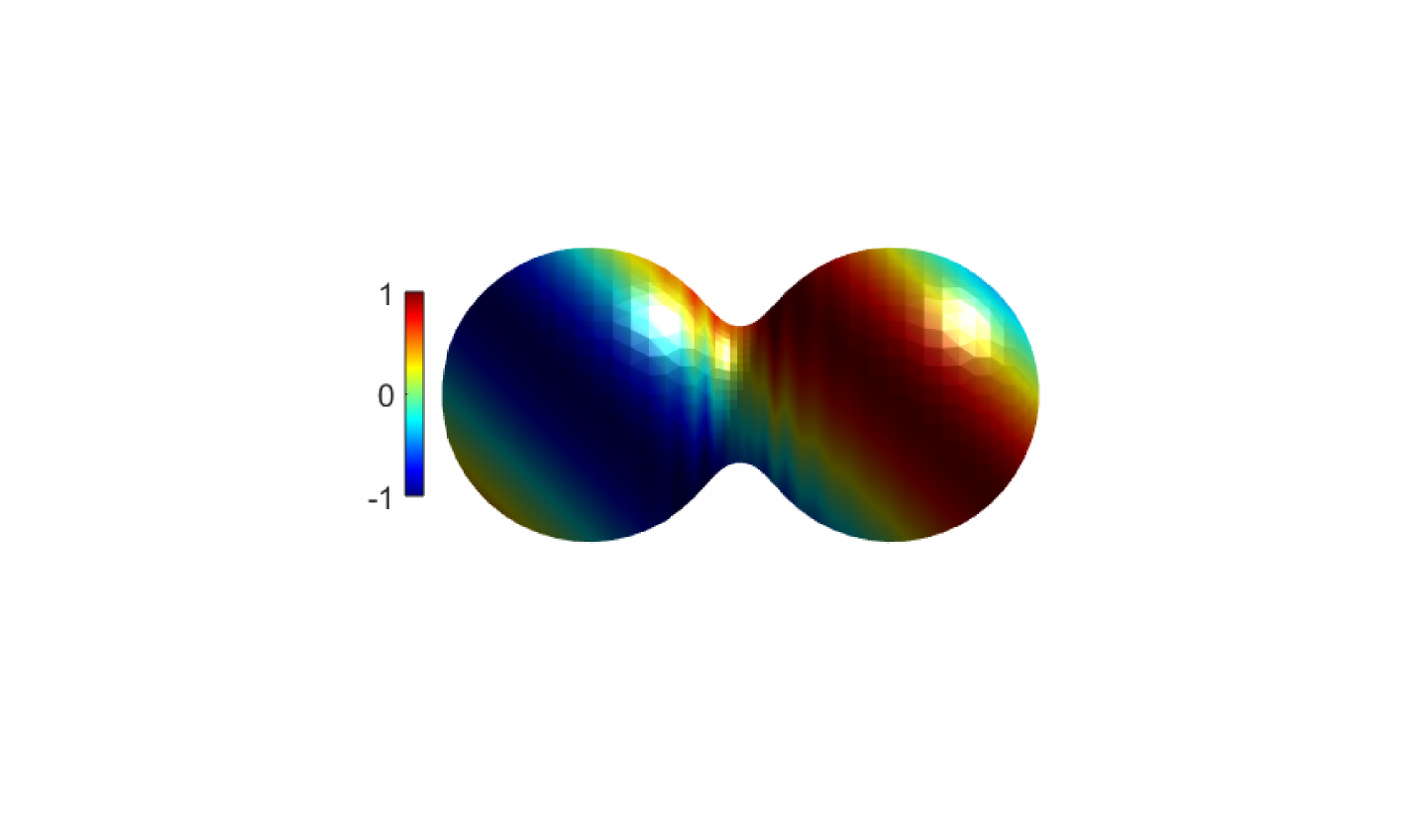}
			\put(43,1.5) {\scriptsize  (a)  $t=0.05$}	
		\end{overpic}
		&	\begin{overpic}[width=0.46\textwidth,trim=100 70 102 60, clip=true,tics=10]{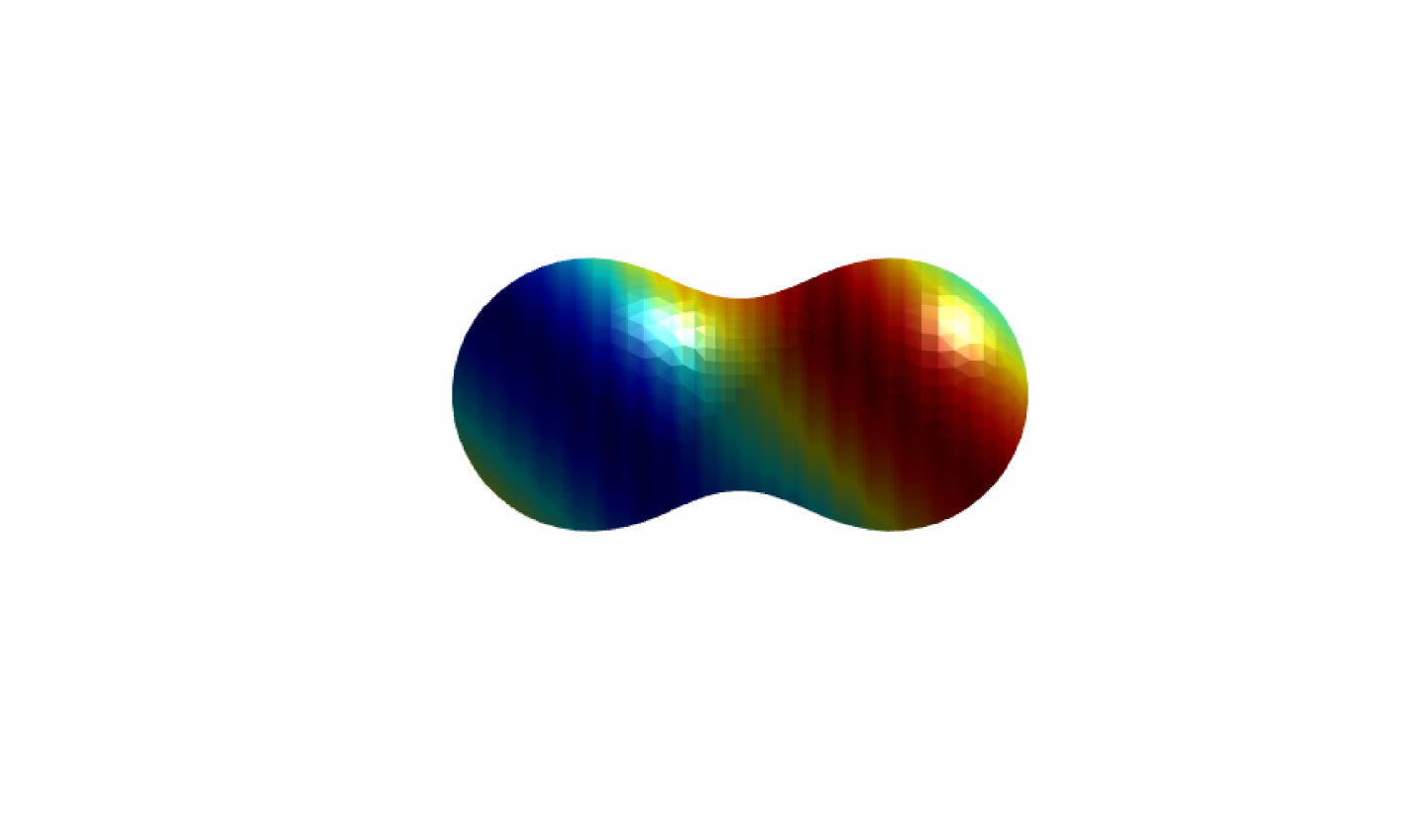}
			\put(43,1.5) {\scriptsize  (b)  $0.2$}
		\end{overpic}
		
	\end{tabular}
	\caption{Example~\ref{Ex_ContSurf}. Numerical solutions at $t=0.05$ and $0.2$, obtained by  Algorithm~\ref{AKM-CN} with exactly determined formulas. We take $\triangle t=0.001$ and the kernels of smoothness order $2$.}\label{Fig_Contin2}
\end{figure}

In  order to test the effect of geometry, we first impose a continuous IC
		$$u^*(\vex,0)=\sin(x+z)$$
	 on the surface. The arrows in  Figure~\ref{Fig_Contin} show mean-curvature velocities and hence the directions of motion.
The diffusion coefficient  and  time step size are set as $\df=1$ and $\triangle t=0.001$ respectively.  First, we solve \eref{eq_CDpupt} with $f=0$ by using the exactly determined formulas with ($n_X=n_Z=1196$).
Figure~\ref{Fig_Contin2} presents two numerical results at $t=0.05$ and $0.2$. From Figure~\ref{Fig_Contin2}(a),  some  oscillations near the region with high curvature can be observed. Once diffusion starts, the oscillatory solution smooth out to the whole surface, see Figure~\ref{Fig_Contin2}(b). We postpone the comparison with oversampling to the next example.

\begin{example}\label{Ex_DiscSurf}{A discontinuous initial condition imposed on point clouds evolving by mean-curvature motion.} \end{example}


\begin{figure}
	\centering	
	\begin{tabular}{cc}
	\multicolumn{2}{c}{	\begin{overpic}[width=0.6\textwidth,trim=100 70 102 60, clip=true,tics=10]{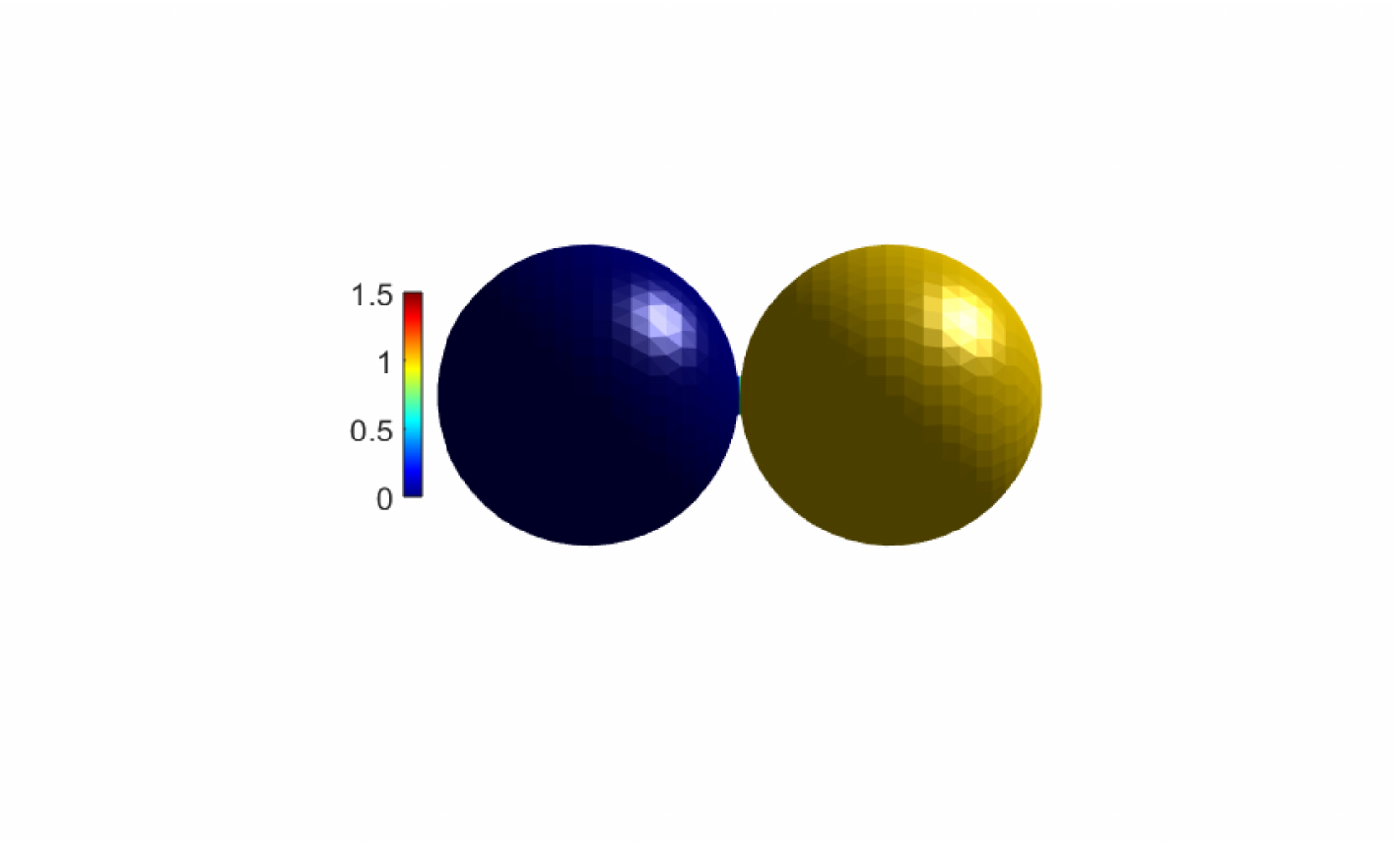}
	\end{overpic}}\\		
	\end{tabular}
    \caption{Example~\ref{Ex_DiscSurf}. Discontinuous initial conditions for a PDE posed on evolving surface under curvature motion.}\label{Fig_Disc4bic}
\end{figure}

We continue the study on the initial high-curvature surface in Example \ref{Ex_ContSurf} but use a discontinuous IC  defined by
\[
    u^*(\vex,0)=0\ (x_1< 0)+1\ (x_1\geq0),
\]
as shown in Figure~\ref{Fig_Disc4bic}.
In the following, we will examine the robustness of Algorithms~\ref{AKM-CN}a and \ref{AKM-CN}b.

\begin{figure}
				\centering	
				\begin{tabular}{cc}
						
				\begin{overpic}[width=0.46\textwidth,trim=100 80 102 50, clip=true,tics=10]{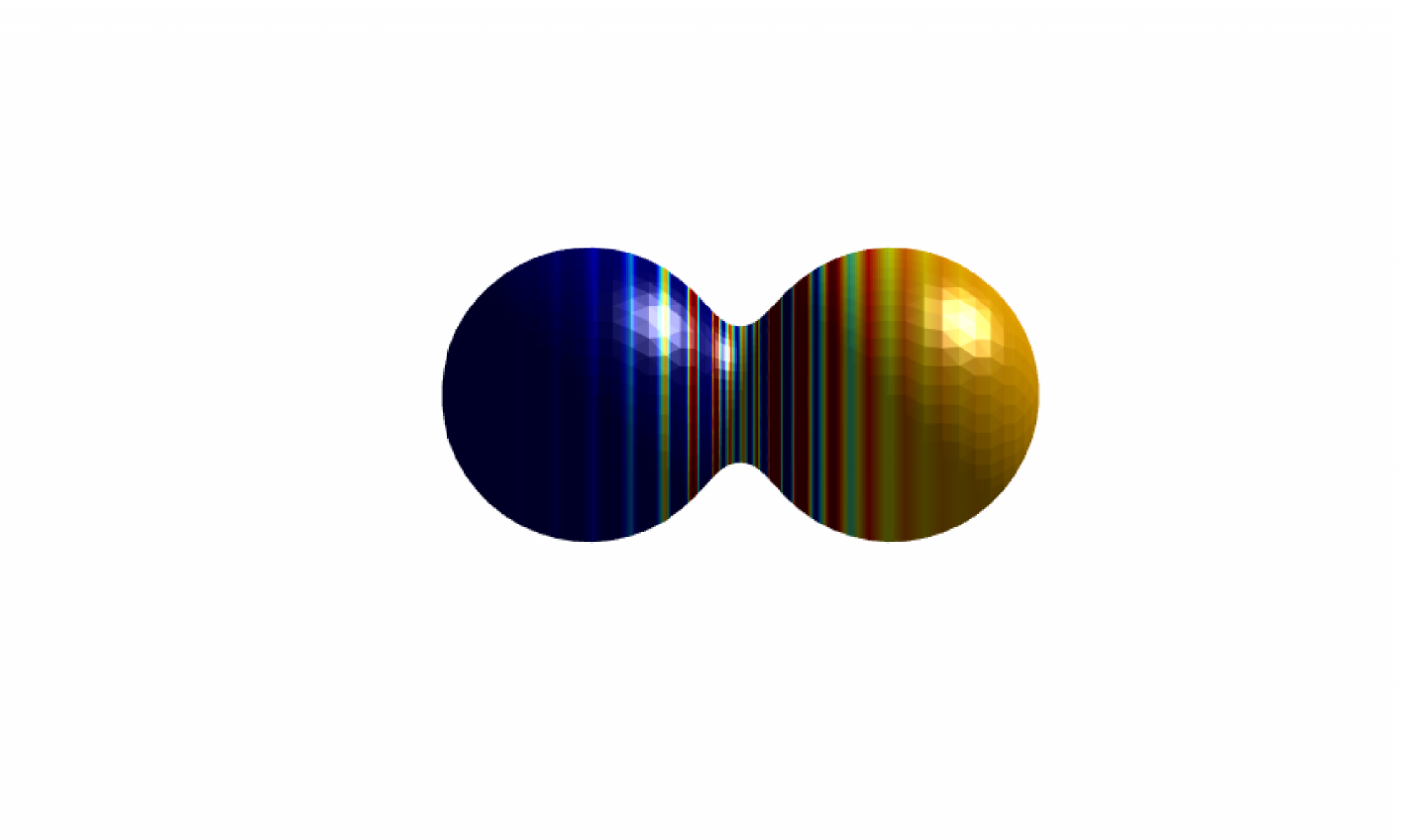}
					\put(45,0) {\scriptsize  (a)  $t=0.05$}	
				\end{overpic}
			&	\begin{overpic}[width=0.46\textwidth,trim=110 150 120 100, clip=true,tics=10]{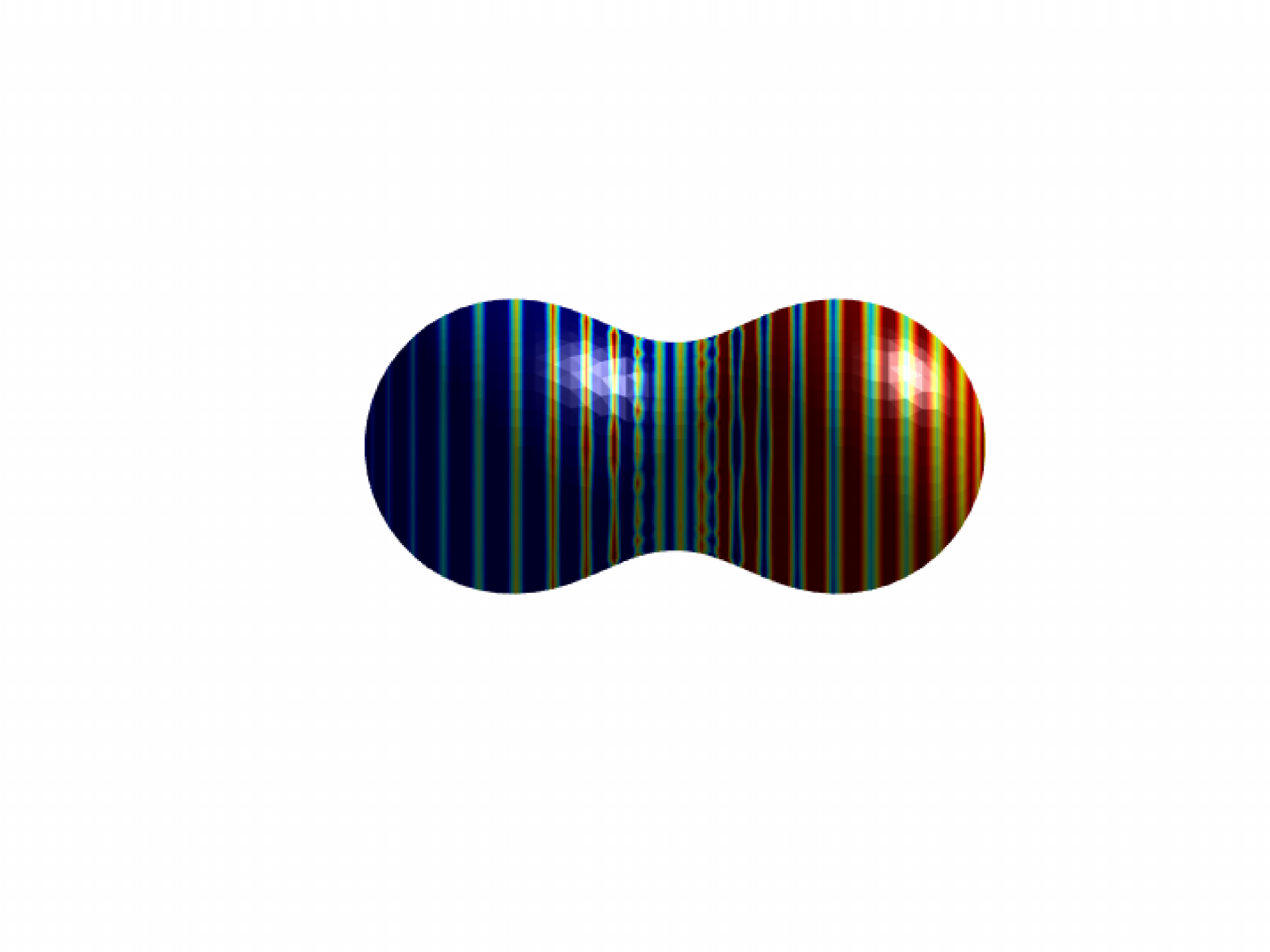}
					\put(50,0) {\scriptsize  (b)  $t=0.2$}
				\end{overpic}\\
			\multicolumn{2}{c}{\scriptsize Exactly determined formula with IC regularization based on $X=Z$ ($n=1196$)}\\
				\begin{overpic}[width=0.46\textwidth,trim=110 150 120 100, clip=true,tics=10]{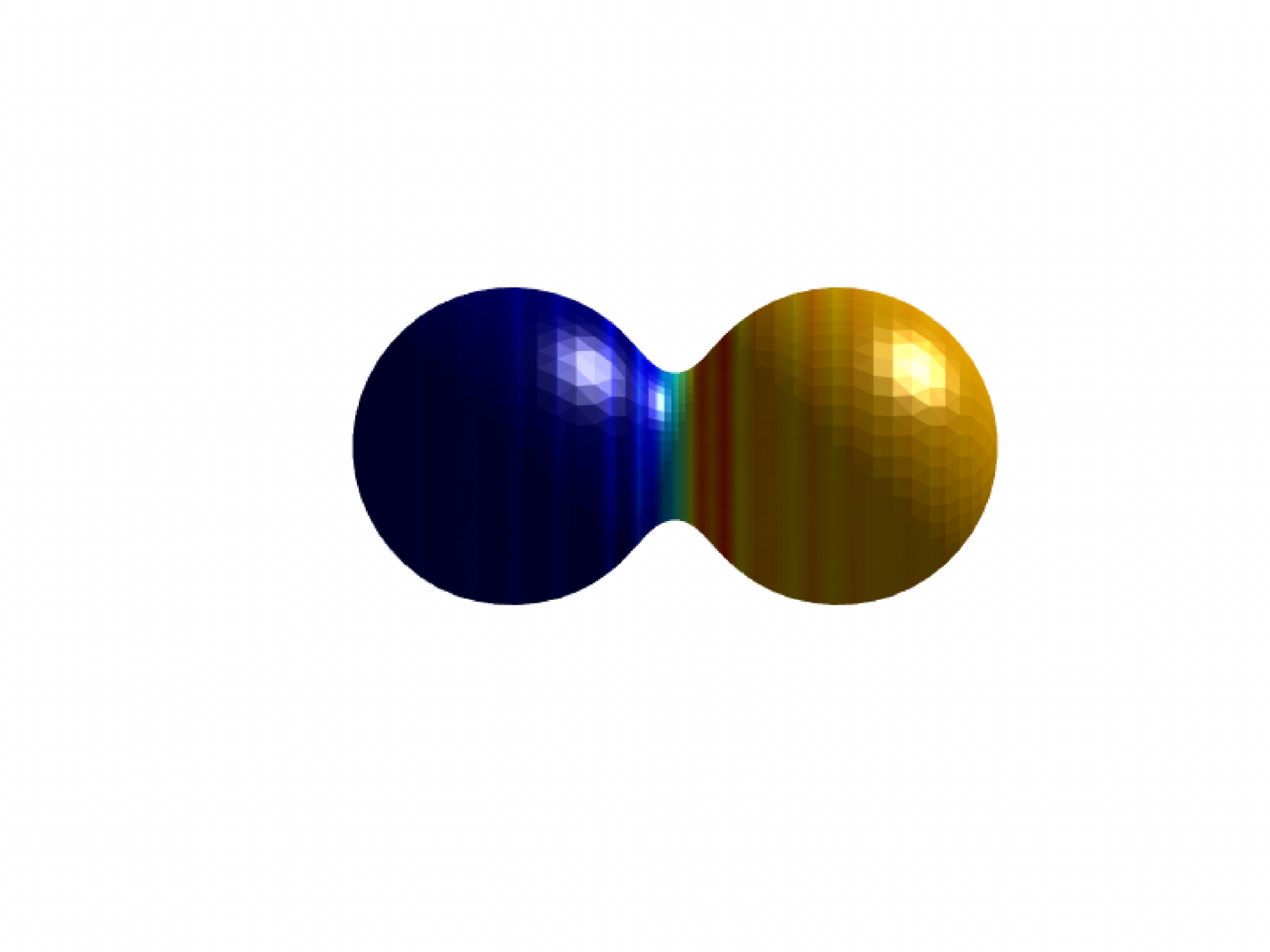}
					
					\put(45,0) {\scriptsize  (c) $t=0.05$}
				\end{overpic}
			&	\begin{overpic}[width=0.46\textwidth,trim=110 150 120 100, clip=true,tics=10]{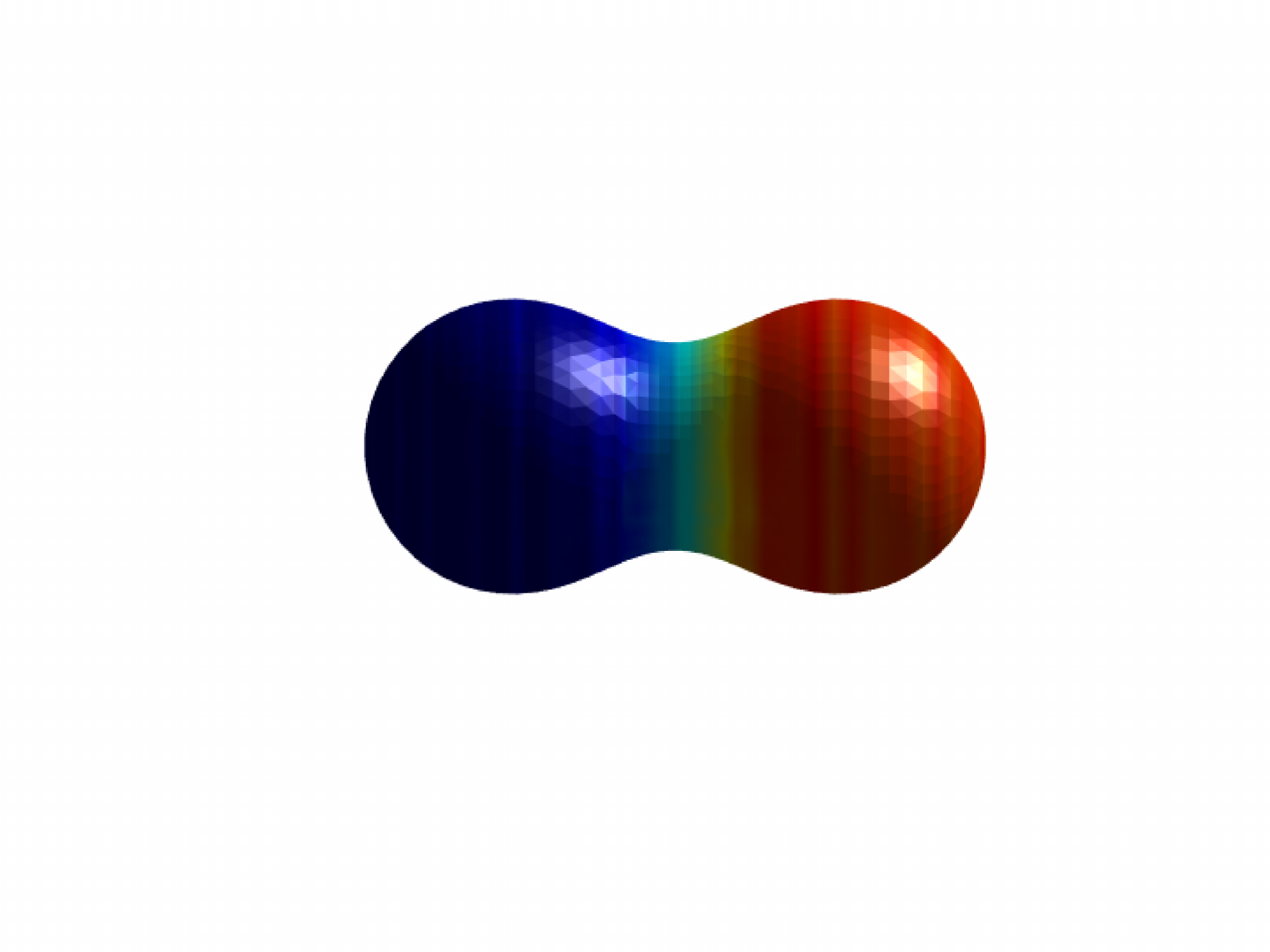}
					\put(50,0) {\scriptsize  (d) $t=0.2$}	\end{overpic}\\
					\multicolumn{2}{c}{\scriptsize Overdetermined formula: Algorithm~\ref{AKM-CN}a with $(n_X,n_Z)=(1736,1196)$}\\
			
					\begin{overpic}[width=0.46\textwidth,trim=110 150 120 100, clip=true,tics=10]{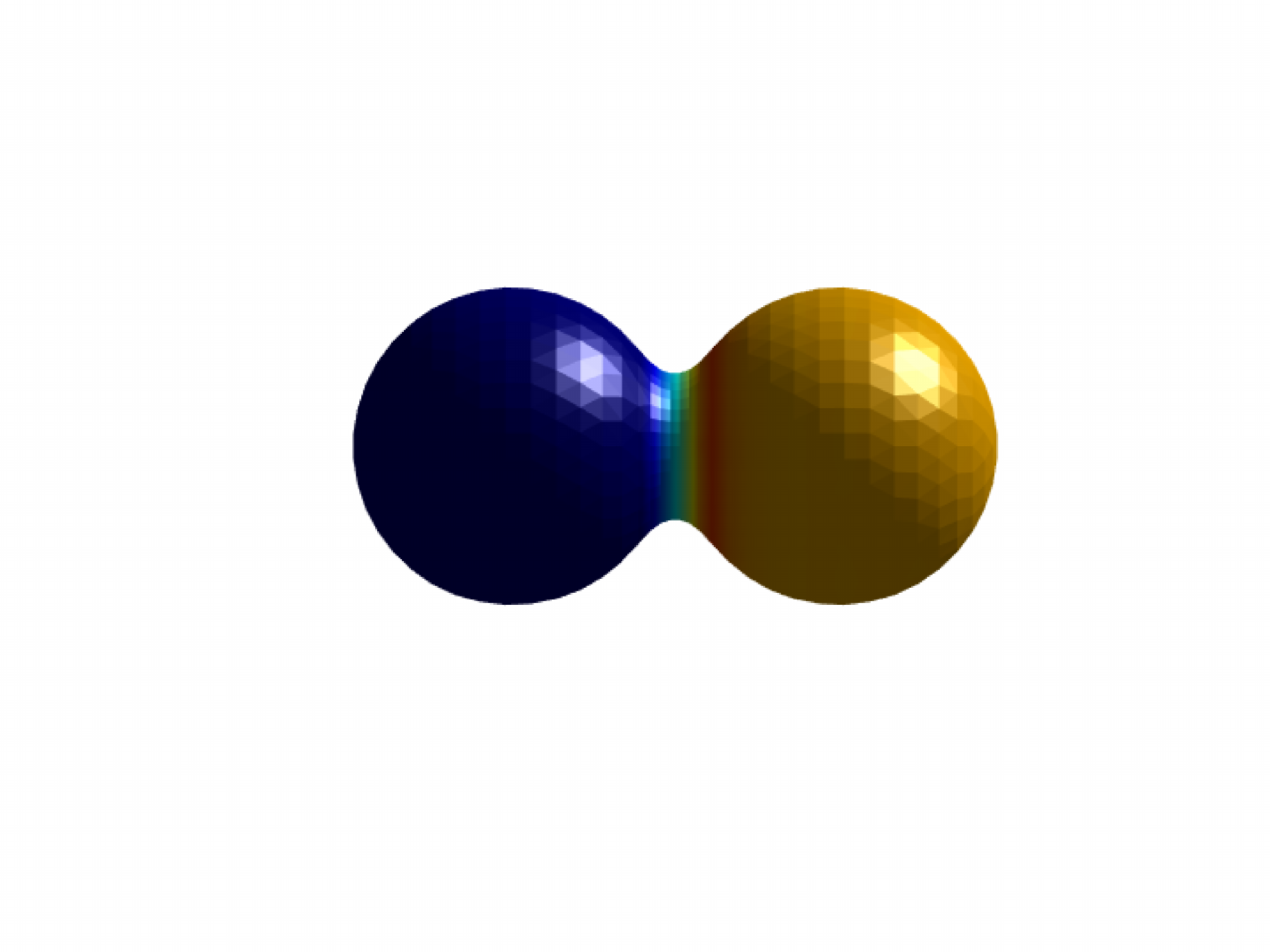}
						
						\put(45,-1) {\scriptsize  (e) $t=0.05$}
					\end{overpic}
					&	\begin{overpic}[width=0.46\textwidth,trim=110 150 120 100, clip=true,tics=10]{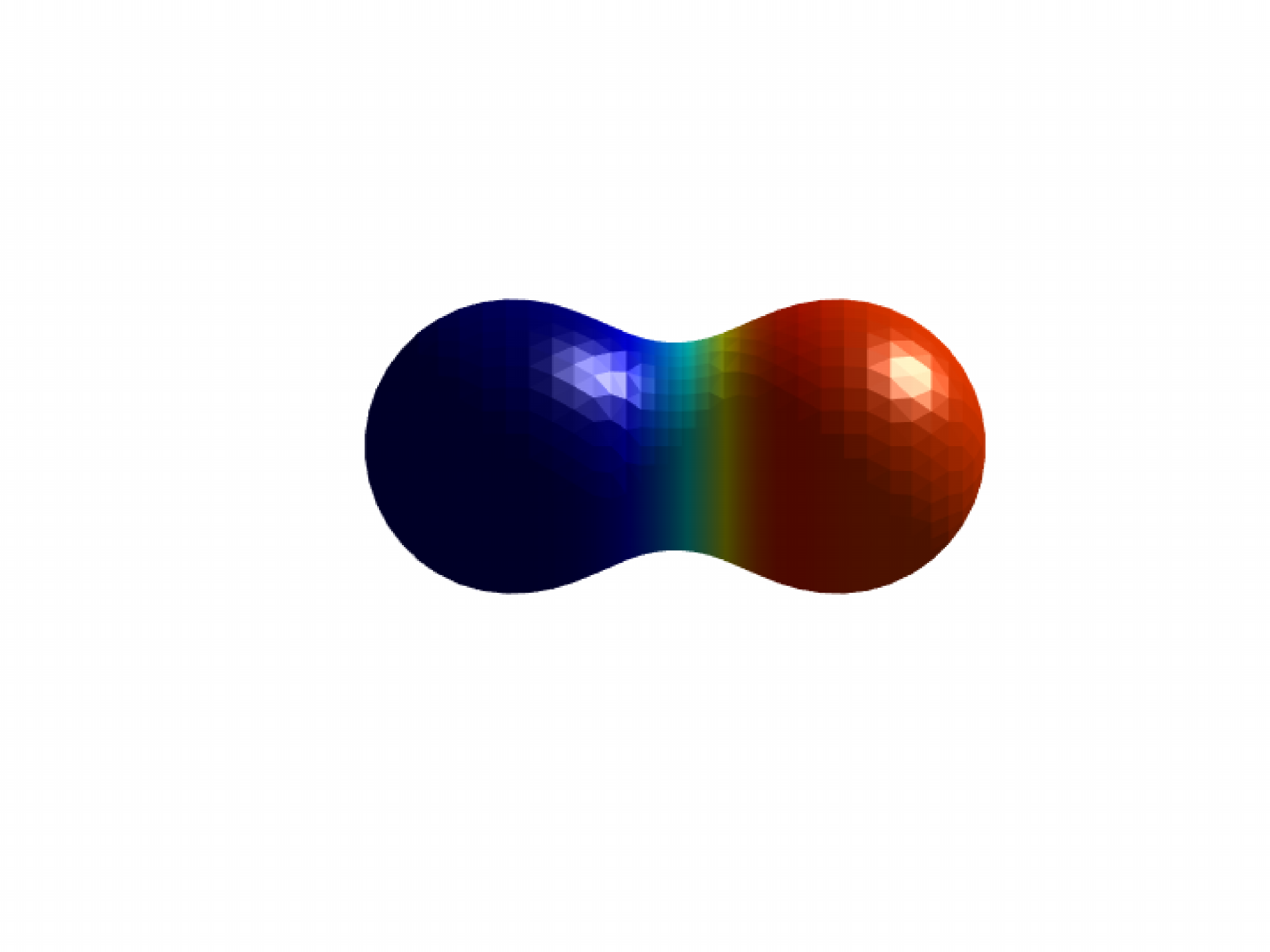}
						\put(50,-1) {\scriptsize  (f) $t=0.2$}
					
				\end{overpic}\\

				\begin{overpic}[width=0.46\textwidth,trim=110 80 120 60, clip=true,tics=10]{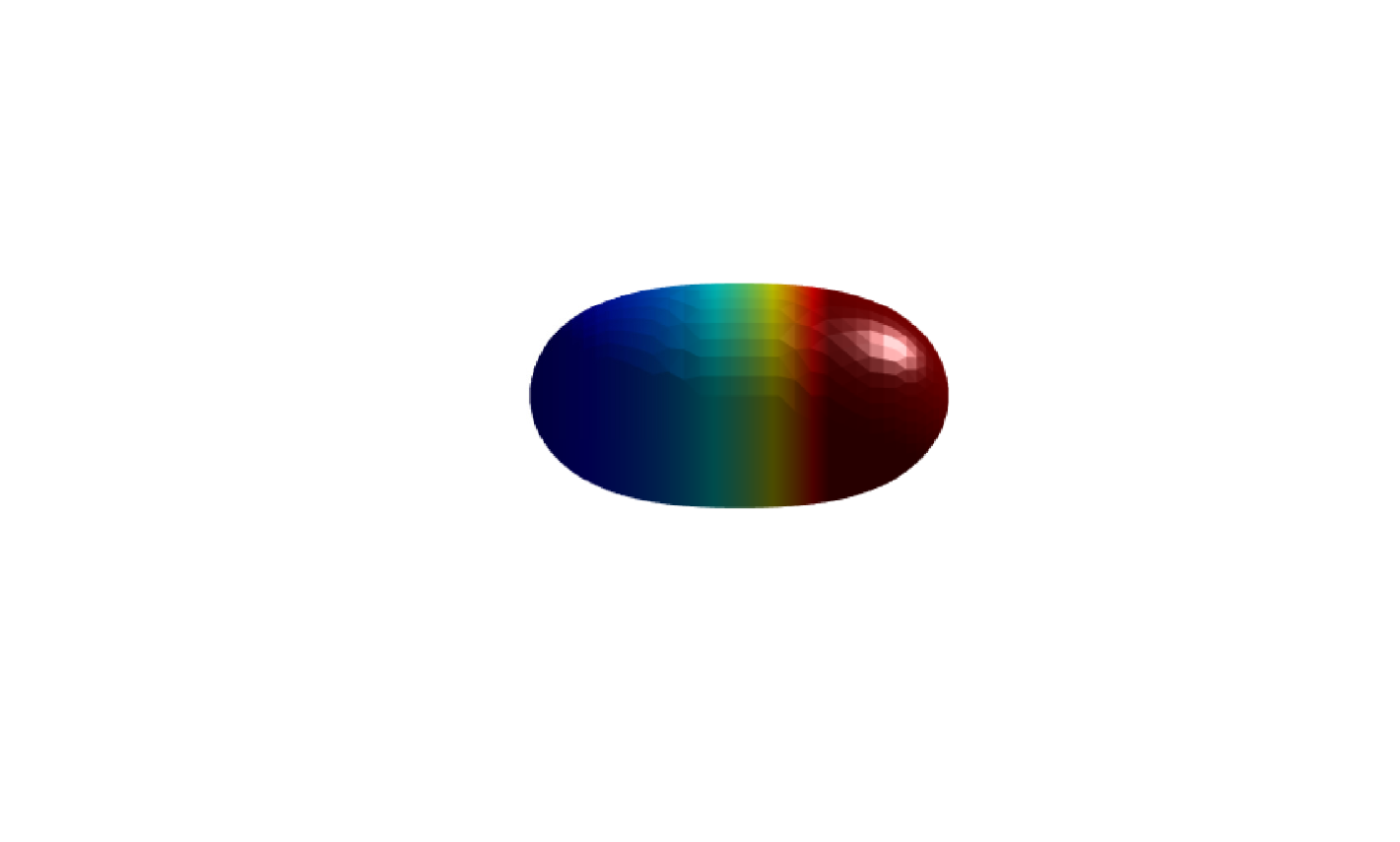}
					\put(48,1.5) {\scriptsize  (g) $t=1$}
				\end{overpic}
		&	\begin{overpic}[width=0.46\textwidth,trim=110 130 120 120, clip=true,tics=10]{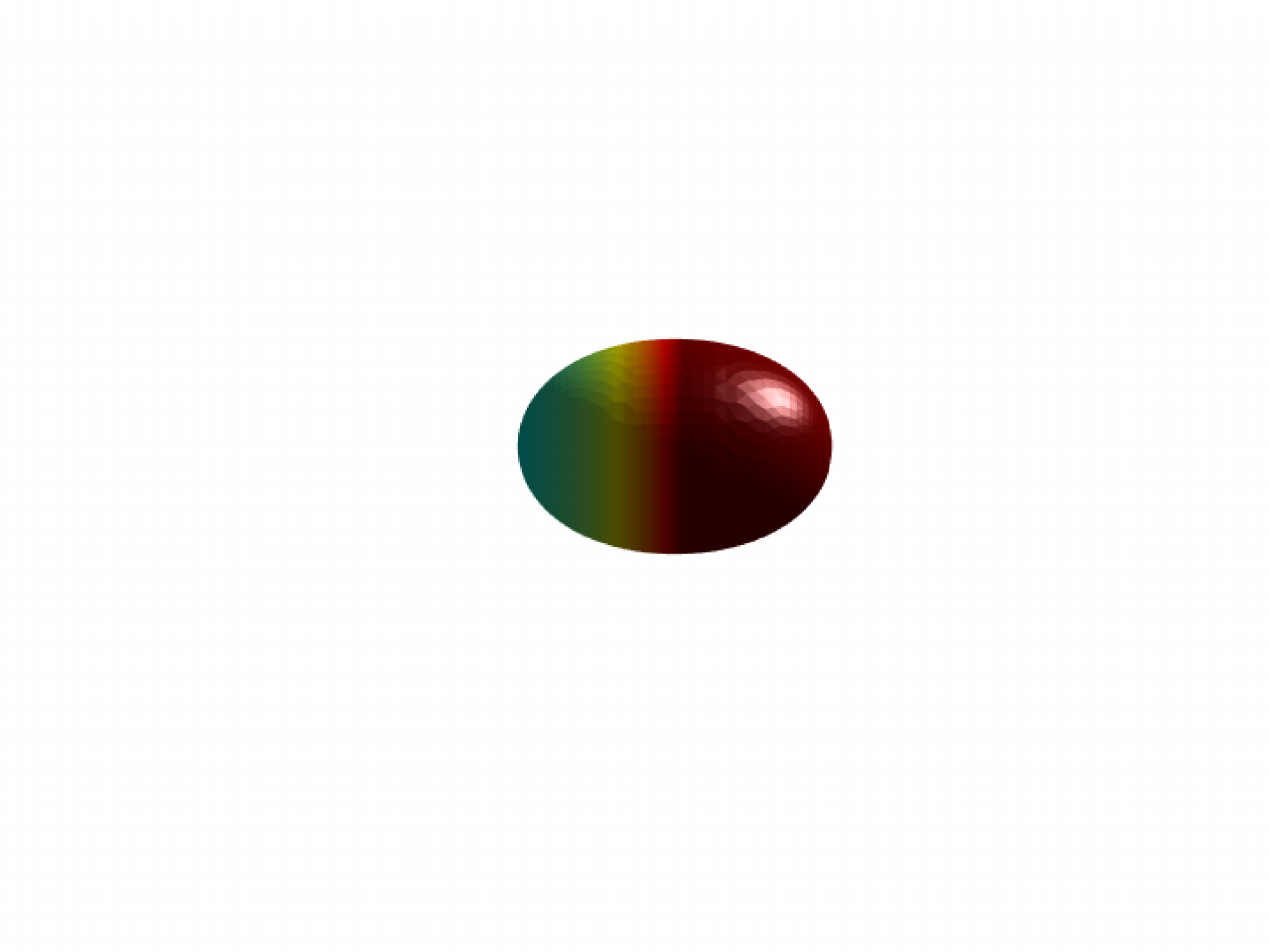}
					\put(50,1.5) {\scriptsize  (h)  $t=1.5$}
				\end{overpic}\\
			\multicolumn{2}{c}{\scriptsize Overdetermined formula: Algorithm~\ref{AKM-CN}b with $(n_X,n_Z)=(1080,710)$}	
			\end{tabular}
				\caption{Example~\ref{Ex_DiscSurf}. Numerical solutions at various $t$, obtained by the exactly determined formulations with IC regularization with a regularization parameter $p= 0.2488$ in  (a)-(b), our Algorithm~\ref{AKM-CN}a in (c)-(d) and Algorithm~\ref{AKM-CN}b in (e)-(h) with different point sets for $X$ and $Z$, under the same settings as in Figure~\ref{Fig_Contin2}.}\label{Fig_Disc4b}
			\end{figure}

The main results were shown in Figure~\ref{Fig_Disc4b}.
Because of the discontinuity,
it makes sense to include some common regularization techniques into our comparison.
We apply a regularized RBF interpolant \cite{fasshauer2007meshfree} to interpolate the IC; then, we update it in time with $X=Z$. The employed regularization parameter $p=0.2488$ is chosen by L-curves in $L^2(\M(t))$ errors \cite{hansen1999curve}.
From Figure~\ref{Fig_Disc4b}(a)-(b), this common-sense setup yields undesirable solutions and we clear see the oscillatory solutions propagating towards the whole surface over time.

Now we focus on the results from Algorithms~\ref{AKM-CN}a and b shown in Figures~\ref{Fig_Disc4b}(c)-(d) and (e)-(h) respectively. As
Algorithm~\ref{AKM-CN}b is a computationally more expensive method, we use it with fewer $Z$ nodes to compare fairly.
To sum up, Algorithm~\ref{AKM-CN}a still show the spurious oscillations initially that will be corrected and the solutions get smoother with increasing time.
Solutions of Algorithm~\ref{AKM-CN}b, in contrast, are physically correct without oscillations.

%

\begin{figure}
				\centering	
				\begin{overpic}[width=0.75\textwidth,trim=60  35  20 15, clip=true,tics=10]{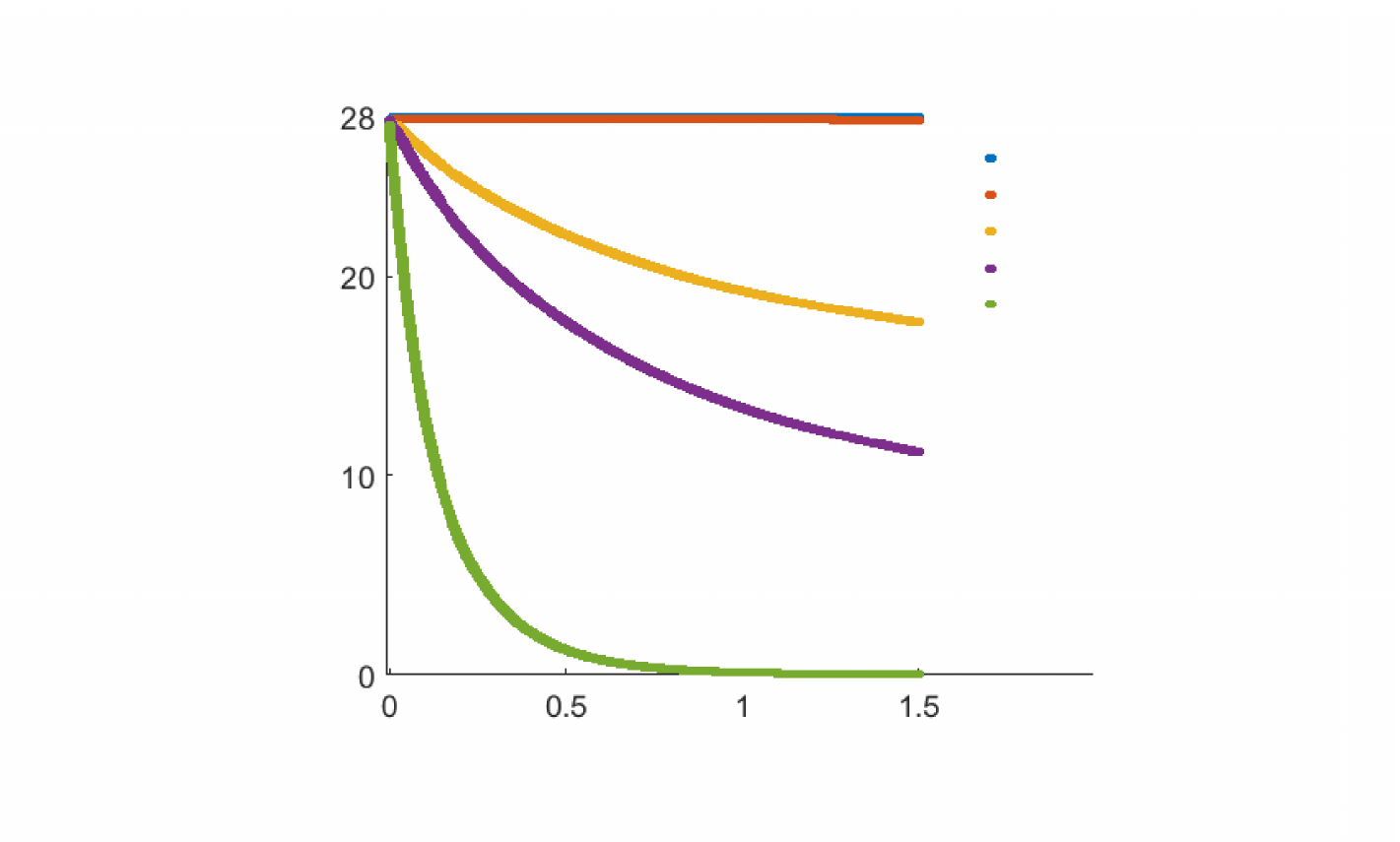}
				\put(6,20) {\scriptsize \rotatebox{90}{Solution of mass}}
				\put(72,49) {\scriptsize Algorithm~\ref{AKM-CN}b}
				\put(72,46) {\scriptsize -- \ref{AKM-CN}a}
				\put(72,42.5) {\scriptsize  Regularized formula with $p=0.2488$}
				\put(72,39.5) {\scriptsize  -- $p=0.3524$}
				\put(72,36.5) {\scriptsize  -- $p=1$}
				\put(48,-3) {\scriptsize  $t$}
				\end{overpic}
				\caption{Example~\ref{Ex_DiscSurf}. Numerical solutions of mass  from $t=0$ to $1.5$, obtained by  Algorithm~\ref{AKM-CN}b and Algorithm~\ref{AKM-CN}a with $(n_X,n_Z)=(1736,1196)$, as well as the regularized formula with $n=1196$ using three different regularization parameters $p=0.2488,\ 0.3524$, and $ 1$ at all times, under the same settings as in Figure~\ref{Fig_Disc4b}.}\label{Fig_Mass_Comp}
			\end{figure}

Before example ends, we want to make sure that Algorithm~\ref{AKM-CN} really is honestly more suitable than regularization.
We now apply regularization at each time, instead of to the initial time only, hoping to avoid oscillatory solutions. Various regularization parameters $p$ were used, all of which yield non-oscillatory smooth \emph{wrong} solutions. We make this conclusion based on mass conservation in Figure~\ref{Fig_Mass_Comp}.

		
			
%

	\section{Conclusions}\label{sec;colclusion}
	
We propose some numerical algorithms  based on intrinsic kernel-based meshless techniques for solving parabolic PDEs on  evolving surfaces.
The proposed algorithms can be implemented either by an analytic projection or a pseudospectral approximation for differential operators on surfaces.
Extension into embedding spaces is unnecessary.
%

The analytic approach (Algorithm~\ref{KPM-CN}) can achieve high-order accuracy and convergence in space when  normal vectors of the surface can be differentiated analytically.
In case of otherwise, the approximated meshless method (Algorithm~\ref{AKM-CN}) becomes handy and it can  handle PDEs defined on point clouds.
Another issue of study is the oversampling strategy, viz., using denser collocation points than trial centers, to yield overdetermined matrix systems to be solved by least-squares. For smooth problems (including initial conditions, surfaces, and solutions), the need of oversampling is not obvious. We demonstrate that, when the initial condition is discontinuous and when the (domain)  surface has high curvature, oversampling is essential for getting physically correct solutions.

%
	
	\bibliographystyle{elsarticle-num}
	\def\cprime{$'$}

\end{document}